\documentclass[epsfig,amstex,amssymb]{amsart}
\usepackage{epsfig}

\begin{document}

\title[Asymptotics]{
        Pseudoholomorphic strips in Symplectisations I: Asymptotic behavior}
\author{
        Casim Abbas     
}
\address{Casim Abbas\\Department of Mathematics\\Michigan State University\\Wells Hall\\East Lansing, MI 48824\\USA}

\date{\today}        
\footnote{This project was partly supported by a New York University Research Challenge Fund Grant and by an NSF Grant (DMS 0196122)}  
\maketitle
\tableofcontents

%
%
\parindent0ex
\parskip1ex plus0.4ex minus0.2ex

\newtheorem{lemma}{Lemma}[section]
\newtheorem{theorem}[lemma]{Theorem}
\newtheorem{corollary}[lemma]{Corollary}
\newtheorem{proposition}[lemma]{Proposition}
\newtheorem{definition}[lemma]{Definition}
\newtheorem{conjecture}{Conjecture}
\newtheorem{example}[lemma]{Example}
\newtheorem{condition}[lemma]{Condition}
\newcommand{\degree}{\mbox{deg}}
\newcommand{\ind}{\mbox{ind}}
\newcommand{\abs}{|}
\newcommand{\ve}{\varepsilon}
\newcommand{\Id}{\mbox{Id}}
\newcommand{\GL}{\mbox{GL}}
\newcommand{\cal}{\mathcal}
\newcommand{\eps}{\varepsilon}
\newcommand{\To}{\longrightarrow}
\newcommand{\Real}{{\bf{R}}}
\newcommand{\Complex}{{\bf{C}}}
\newcommand{\RM}{{\bf{R}}\times M}
\newcommand{\RL}{{\bf{R}}\times {\mathcal L}}
\newcommand{\tu}{\tilde{u}}
\newcommand{\pil}{\pi_{\lambda}}
\newcommand{\pas}{\partial_s}
\newcommand{\pat}{\partial_t}
\newcommand{\od}{\{0\}\times{\mathcal D}^{\ast}}
\newcommand{\SC}{{\mathcal S}}
\newcommand{\ho}{\mbox{Hom}}
\newcommand{\oo}{{\mathcal O}}

\begin{abstract}
This paper is part of a larger program, the investigation of the Chord Problem in three dimensional contact geometry. The main tool will be pseudoholomorphic strips in the symplectisation of a three dimensional contact manifold with two totally real submanifolds $L_0,L_1$ as boundary conditions. The submanifolds $L_0$ and $L_1$ do not intersect transversally. The subject of this paper is to study the asymptotic behavior of such pseudoholomorphic strips. 
\end{abstract}         
\section{Introduction}
This paper is the first part of a larger program, the investigation of the chord problem in three dimensional contact geometry (\cite{part2}, \cite{part3}, \cite{part4}).
Let $(M,\lambda)$ be a $2n+1$--dimensional contact manifold, i.e. $\lambda$ is a 1--form on $M$ such that $\lambda\wedge (d\lambda)^n$ is a volume form on $M$. The contact structure associated to $\lambda$ is the $2n$--dimensional vector bundle $\xi=\ker\lambda\rightarrow M$, which is a symplectic vector bundle  with symplectic structure $d\lambda|_{\xi\oplus\xi}$. There is a distinguished vector field associated to a contact form, the Reeb vector field $X_{\lambda}$, which is defined by the equations 
\[
i_{X_{\lambda}}d\lambda\equiv 0\ ,\ i_{X_{\lambda}}\lambda\equiv 1.
\]
We denote by $\pil:TM\rightarrow\xi$ the projection along the Reeb vector field.
The Chord Problem is about the global dynamics of the Reeb vector field. More precisely, the issue is the existence of so--called 'characteristic chords'. These are trajectories $x$ of the Reeb vector field which hit a given Legendrian submanifold ${\mathcal L}^n\subset (M,\lambda)$ at two different times $t=0$ , $T>0$. We also ask for $x(0)\neq x(T)$, otherwise the chord would actually be a periodic orbit. Recall that a submanifold $L$ in a $2n+1$--dimensional contact manifold $(M,\lambda)$ is called Legendrian if it is everywhere tangent to the hyperplane field $\xi$ and if it has dimension $n$. We are mostly interested in the three--dimensional situation, the question is then whether a given Legendrian knot has a characteristic chord. The Chord problem should be viewed as the relative version of the Weinstein conjecture which deals with the existence of periodic orbits of the Reeb vector field.\\
Characteristic chords occur naturally in classical mechanics. In this context they are referred to as 'brake--orbits', and were investigated by Seifert in 1948 \cite{Seifert} and others since the 1970's \cite{Ambrosetti}, \cite{Bolotin}, \cite{VanGroesen}, \cite{Weinstein}.\\
In 1986, V.I. Arnold conjectured the existence of characteristic chords on the three sphere for any contact form inducing the standard contact structure and for any Legendrian knot \cite{Arnold}. After a partial result by the author in \cite{A3} this conjecture was finally confirmed by K. Mohnke in \cite{Mohnke}. It is natural to ask the existence question for characteristic chords not only for $M=S^3$, but also for general contact manifolds. A new invariant for Legendrian knots and contact manifolds proposed by Y. Eliashberg, A. Givental and H. Hofer in \cite{SFT} ('Relative Contact Homology') is actually based on counting characteristic chords and periodic orbits of the Reeb vector field.\\
The subject of the paper \cite{part4} is an existence result for characteristic chords which goes beyond the special classes of contact three manifolds investigated so far. The purpose of this paper and \cite{part2}, \cite{part3} is to establish a filling method by pseudoholomorphic curves where we use a surface $F\subset M=M^3$ with boundary, and where we start filling from a tangency at the boundary. Pseudoholomorphic curves are maps from a Riemann surface into an almost complex manifold $W$ satisfying a nonlinear Cauchy Riemann type equation. In our case, the manifold $W$ is the symplectisation $({\bf R}\times M,d(e^t)\lambda)$ of the contact manifold $(M,\lambda)$. We are going to consider a special type of almost complex structures $\tilde{J}$ on ${\bf R}\times M$. We pick a complex structure $J:\xi\rightarrow\xi$ such that $d\lambda\circ(\mbox{Id}\times J)$ is a bundle metric on $\xi$. We then define an almost complex structure on ${\bf R}\times M$ by demanding $\tilde{J}\equiv J$ on $\xi$ and sending $\partial/\partial t$ (the generator of the ${\bf R}$--component) onto the Reeb vector field. Then $\tilde{J}(p)$ has to map $X_{\lambda}(p)$ onto $-\partial/\partial t$.\\
If $S$ is a Riemann surface with complex structure $j$ then we define a 
map 
$$
\tilde{u}=(a,u):S\longrightarrow{\bf R}\times M
$$
to be a pseudoholomorphic curve if
$$
D\tilde{u}(z)\circ j(z)=\tilde{J}(\tilde{u}(z))\circ D\tilde{u}(z)\ \mbox{for all}\ z\in S.
$$
If $(s,t)$ are conformal coordinates on $S$ then this becomes:
\[
\partial_s\tilde{u}+\tilde{J}(\tilde{u})\partial_t\tilde{u}=0.
\]
We are interested only in pseudoholomorphic curves which have finite energy in the sense that
$$
E(\tilde{u}):=\sup_{\phi\in\Sigma}\int_{S}
\tilde{u}^{\ast}d(\phi\lambda) <+\infty
$$
where $\Sigma:=\{\phi\in C^{\infty}({\bf R},[0,1])\,|\,\phi'\ge 
0\}$. The Riemann surface $S$ in this paper is an infinite strip $S={\bf R}\times [0,1]$, and we will impose a mixed boundary condition as follows: Let ${\mathcal L}\subset M$ be a homologically trivial Legendrian knot bounding an embedded surface ${\mathcal D}$. A point $p\in{\mathcal D}$ is called singular if $T_p{\mathcal D}=\ker\lambda(p)$. If the surface is oriented (by a volume form $\sigma$) and if $j:{\mathcal D}\hookrightarrow M$ is the inclusion, then we define a vector field $Z$ on ${\mathcal D}$ by $i_Z\sigma=j^{\ast}\lambda$. This vector field vanishes precisely in the singular points. The flow lines of $Z$ determine a singular foliation of the surface ${\mathcal D}$ which does not depend on the particular choice of the volume form or the contact form. This singular foliation is also called the characteristic foliation of ${\mathcal D}$ (induced by $\ker\lambda$). Let $p\in{\mathcal D}$ be a singular point and denote by $Z'(p):T_p{\mathcal D}\rightarrow T_p{\mathcal D}$ the linearization of the vector field $Z$ in $p$. Let $\lambda_1,\lambda_2$ be the eigenvalues of $Z'(p)$. We say that $p$ is non--degenerate if none of the eigenvalues lie on the imaginary axis. A non--degenerate singular point $p$ is called elliptic if $\lambda_1\lambda_2>0$ and hyperbolic if $\lambda_1\lambda_2<0$. In the elliptic case the critical point $Z(p)=0$ is either a source or a sink, and in the hyperbolic case it is a saddle point.\\  
Choosing ${\mathcal D}$ appropriately we may assume that there are only non--degenerate singular points, in particular there are only finitely many. We denote the surface without the singular points by ${\mathcal D}^{\ast}$. We consider the boundary value problem 
\begin{equation}\label{main-boundary-value-problem}
\begin{array}{ll}
          \tilde{u}=(a,u):S\longrightarrow\RM & \\
          \partial_s\tu+\tilde{J}(\tu)\partial_t\tu=0 & \\
          \tu(s,0)\subset{\bf R}\times{\cal L} & \\
          \tu(s,1)\subset\{0\}\times{\cal D}^{\ast} & \\
           0< E(\tu)<+\infty. & 
          \end{array}
\end{equation}
The subject of this paper is to investigate the behavior of solutions $\tu$ for large $|s|$. The finiteness condition on the energy actually forces the solutions to converge to points $\tilde{p}_{\pm}\in\{0\}\times{\mathcal L}$ at an exponential rate. We first introduce suitable coordinates near the Legendrian knot, and we deform ${\mathcal D}$ near its boundary, keeping ${\mathcal L}=\partial {\mathcal D}$ fixed, in order to achieve a certain normal form for ${\mathcal D}$ near its boundary. We then derive exponential decay estimates for $\tilde{u}-\tilde{p}_{\pm}$ and all its derivatives in coordinates. In local coordinates near $\tilde{p}_{\pm}$ the almost complex structure $\tilde{J}$ corresponds to some real $4\times 4$--matrix valued function which we denote by $M$. The main result of this paper is the following asymptotic formula
\begin{theorem}\label{asymptotic-formula-theorem}
For sufficiently large $s_0$ and $s\ge s_0$ we have the following
asymptotic formula for non constant solutions $v$ of
(\ref{main-boundary-value-problem}) having finite energy:
\begin{equation}
v(s,t)=e^{\int_{s_0}^s\alpha(\tau)d\tau}\Big(e(t)+r(s,t)\Big),
\end{equation}
where $\alpha:[s_0,\infty)\rightarrow{\bf R}$ is a smooth function
satisfying $\alpha(s)\rightarrow\lambda < 0$ as
$s\rightarrow\infty$ with $\lambda$ being an eigenvalue of the
selfadjoint operator
\[
A_{\infty}:L^2([0,1],{\bf R}^4)\supset H^{1,2}_L([0,1],{\bf
R}^4)\longrightarrow L^2([0,1],{\bf R}^4)
\]
\[
\gamma\longmapsto -M_{\infty}\dot{\gamma}\ ,\
M_{\infty}:=\lim_{s\rightarrow\infty}M(v(s,t)).
\]
Moreover, $e(t)$ is an eigenvector of $A_{\infty}$ belonging to
the eigenvalue $\lambda$ with $e(t)\neq 0$ for all $t\in[0,1]$, 
and $r$ is a smooth function so that $r$
and all its derivatives converge to zero uniformly in $t$ as
$s\rightarrow\infty$. 
\end{theorem}
We will prove more about the decay of $|\lambda-\alpha(s)|$, $r$ and their derivatives:
\begin{theorem}\label{convergence-of-alpha}
Let $r$ and $\alpha(s)$ be as in theorem \ref{asymptotic-formula-theorem}. Then there is a constant $\delta>0$ such that for each integer $l\ge 0$ and each multi--index $\beta\in{\bf N}^2$
\[
\sup_{0\le t\le 1}|D^{\beta}r(s,t)|\ ,\ \left|\frac{d^l}{ds^l}(\alpha(s)-\lambda)\right|\le c_{\beta,l} \,e^{-\delta|s|}
\]
with suitable constants $c_{\beta,l}>0$.
\end{theorem}
The subscript 'L' in $H^{1,2}_L([0,1],{\bf R}^4)$ indicates the boundary condition (see (\ref{definition-of-H12L}) for a precise definition). This formula is an essential ingredient for the rest of the program \cite{part2}, \cite{part3}, \cite{part4}.\\
The asymptotic behavior of holomorphic strips with mixed boundary conditions similar to ours was investigated in \cite{RS2}, but only for non--degenerate ends. We are dealing with a degenerate situation, i.e. the manifolds $L_0={\bf R}\times {\mathcal L}$ and $L_1=\{0\}\times {\mathcal D}$ do not intersect transversally. The degenerate situation is much more delicate: In the non--degenerate case the intersection $L_0\cap L_1$ would consist of isolated points. Having shown that a pseudoholomorphic strip $\tu(s,t)$ with finite energy approaches $L_0\cap L_1$ as $|s|\rightarrow\infty$ one can fairly easily see that oscillations between two points in $L_0\cap L_1$ would cost too much energy, i.e. it would contradict $E(\tu)<\infty$. In our case we have to show that the end of the solution cannot move along the 1--dimensional set $L_0\cap L_1$ while $|s|$ grows. Analytically, degeneracy means that the operator $A_{\infty}$ above has a nontrivial kernel. The strategy is to derive estimates for the 'components' of $\tu$ orthogonal to the kernel of $A_{\infty}$ (in a suitable sense). We will then show that they decay fast enough to force the component along the kernel of $A_{\infty}$ to zero as well.\\  
Degenerate ends were investigated in the paper \cite{HWZIV}, but only for pseudoholomorphic cylinders $S={\bf R}\times S^1$ (periodic boundary condition in $t$). Our problem requires a different approach. The paper \cite{RS2} contains the decay estimate of theorem \ref{convergence-of-alpha} for the case $\beta=0$. Eduardo Mora proved theorem \ref{convergence-of-alpha} for pseudoholomorphic cylinders independently of the author in his Ph.D. thesis \cite{Mora}.\\
Because we are choosing special $\tilde{J}$ and ${\mathcal D}$ near $\{0\}\times{\mathcal L}$, we obtain as a side product the following existence statement
\begin{theorem}\label{local-existence-theorem}
Let $(M,\lambda)$ be a three dimensional contact manifold. Moreover, let ${\mathcal L}$ be a Legendrian knot which bounds an embedded surface ${\mathcal D}'$ so that the characteristic foliation has only finitely many singular points. Then there is another embedded surface ${\mathcal D}$ which is a $C^0$ small perturbation of ${\mathcal D}'$ having the same boundary and the same singular points as ${\mathcal D}'$ and a $d\lambda$--compatible complex structure $J:\ker\lambda\rightarrow\ker\lambda$ so that the following is true: Near each elliptic singular point $e\in\partial{\mathcal D}={\mathcal L}$ there are embedded solutions $\tu_{\tau}$ , $0<\tau<1$ to the boundary value problem (\ref{main-boundary-value-problem}) with the properties:
\begin{itemize}
\item $\tu_{\tau}(S)\cap\tu_{\tau'}(S)=\emptyset$ if $\tau\neq\tau'$,
\item $\tu_{\tau}\rightarrow e$ uniformly with all derivatives as $\tau\rightarrow 0$,
\item the family $\tu_{\tau}$ depends smoothly on the parameter $\tau$.
\end{itemize}
\end{theorem}
\section{Simplifying the spanning surface ${\mathcal D}$ near the
boundary}\label{normal-form-for-the-surface} In this section we will
simplify the surface ${\mathcal D}$ near its boundary to obtain a
normal form in coordinates near the knot ${\mathcal L}=\partial{\mathcal D}$. This is useful for the analysis later. In particular, we will be able to produce explicitly a family of finite energy strips coming out from elliptic singular points on the boundary.\\ 
If $(M,\lambda)$ is a three dimensional contact manifold and ${\mathcal L}$ a Legendrian knot
in $M$ then, by a well--known theorem of A. Weinstein (see
\cite{Weinstein2}, \cite{Weinstein3} and \cite{AH}), there are open neighborhoods
$U\subset M$ of the knot ${\mathcal L}$, $V\subset S^1\times
\Real^2$ of $S^1\times\{(0,0)\}$ and a diffeomorphism
$\Psi:U\rightarrow V$, so that
$\Psi^{\ast}(dy+xd\theta)=\lambda|_U$, where $\theta$ denotes the
coordinate on $S^1\approx\Real/{\bf{Z}}$ and $x,y$ are coordinates
on $\Real^2$. We will refer to this result as the 'Legendrian neighborhood theorem'. If we are working near the knot ${\mathcal L}$ we
may assume that our contact manifold is
$(S^1\times\Real^2,\lambda=dy+xd\theta)$ and the knot is given by
$S^1\times\{(0,0)\}$. We will denote the piece of the spanning
surface ${\mathcal D}\cap U$ again by ${\mathcal D}$. Choosing $U$ sufficiently small we may assume that all the singular points on the piece ${\mathcal D}\cap U$ lie on the boundary, i.e.
\[
\{p\in {\mathcal D}\cap U\,|\,T_p{\mathcal D}=\ker\lambda(p)\}=\{(\theta_k,0,0)\in S^1\times {\bf R}^2\}_{1\le k\le N}\ ,\ N\in{\bf N}.
\]
We parameterize ${\mathcal D}$ as
follows:
\[
{\mathcal D}=\{(\theta,x(\theta,r),y(\theta,r))\in
S^1\times\Real^2\,|\,r,\theta\in[0,1]\},
\]
where $x,y$ are suitable smooth functions which are 1--periodic in
$\theta$ and satisfy
\[
x(\theta,0)\equiv y(\theta,0)\equiv 0.
\]
Moreover we orient ${\mathcal D}$ in such a way that the above
parameterization $([0,1]\times[0,1],d\theta\wedge
dr)\rightarrow{\mathcal D}$ is orientation preserving. We orient
${\mathcal L}$ by $v=\frac{d}{d\theta}$, so that $(v,\nu)$ is
positively oriented, where $\nu$ denotes the inward normal vector.
A point $(\theta_0,0,0)$ is a singular point if and only if
$\partial_ry(\theta_0,0)=0$. Since also
$\partial_{\theta}y(\theta_0,0)=\partial_{\theta}x(\theta_0,0)=0$
and ${\mathcal D}$ is embedded, we conclude that
$\partial_rx(\theta_0,0)\neq 0$. The tangent space
$T_{(\theta_0,0,0)}{\mathcal D}$ is oriented by the basis
$\left(\frac{\partial}{\partial\theta},\partial_rx(\theta_0,0)\frac{\partial}{\partial
x}\right)$. On the other hand, the contact structure
$\ker\lambda(\theta_0,0,0)$ is oriented by
$\left(\frac{\partial}{\partial\theta},-\frac{\partial}{\partial
x}\right)$. The singular point $(\theta_0,0,0)$ is called positive
if these two orientations coincide, which is the case for
$\partial_rx(\theta_0,0)<0$, otherwise $(\theta_0,0,0)$ is called
negative. Hence in the case of a positive (negative) singularity,
the surface ${\mathcal D}$ lies on the side of the negative
(positive) x--axis. We would like to perturb ${\mathcal D}$ near
its boundary, leaving the boundary fixed, so that the number and
type of the singularities does not change and the new surface has
some kind of normal form near its boundary. the following is the
main result of this section. It is an immediate consequence of
propositions \ref{3.1.2.} and \ref{3.1.3.} below.\\

\begin{proposition}\label{3.1.1.}
Let $(M,\lambda)$ be a three-dimensional contact manifold. Further, let ${\mathcal L}$ be
a Legendrian knot and ${\mathcal D}$ an embedded surface with
$\partial{\mathcal D}={\mathcal L}$ so that all the singular points are non--degenerate. We denote the finitely many singular points on the boundary by $e_k$ , $1\le k\le N$ (ordered by moving in the direction of the orientation of ${\mathcal L}$).\\ 
Then there is an embedded surface ${\mathcal D}'$
having the same boundary as ${\mathcal D}$ which differs from
${\mathcal D}$ only by a $C^0$--small perturbation supported near
${\mathcal L}$ having the same singular points as ${\mathcal D}$ so that the following holds:\\
There is a neighborhood $U$ of ${\mathcal L}$ and a
diffeomorphism $\Phi:U\rightarrow S^1\times{\bf R}^2$ so that
\begin{itemize}
\item $\Phi^{\ast}(dy+xd\theta)=\lambda|_U\ ,\ (\theta,x,y)\in
S^1\times{\bf R}^2,$
\item $\Phi({\mathcal L})=S^1\times\{(0,0)\},$
\item $\Phi(e_k)=(\theta_k,0,0)$\ ,\ $0\le\theta_1<\cdots<\theta_N<1$,
\item $\Phi(U\cap{\mathcal
D}')=\{(\theta,a(\theta)r,b(\theta)r)\in S^1\times{\bf
R}^2\,|\,\theta,r\in[0,1]\},$
\end{itemize}
where $a,b$ are smooth 1--periodic functions with:
\begin{itemize}
\item $b(\theta_k)=0$ and $b(\theta)$ is nonzero if $\theta\neq\theta_k$,
\item $a(\theta_k)<0$ if $e_k$ is a positive singular point, $a(\theta_k)>0$ if $e_k$ is a negative singular point,
\item if $e_k$ is elliptic then $-1<\frac{b'(\theta_k)}{a(\theta_k)}<0$,
\item if $e_k$ is hyperbolic then the quotient $\frac{b'(\theta_k)}{a(\theta_k)}$ is either strictly smaller than $-1$ or positive,
\item $a$ has exactly one zero in each of the intervals
$[\theta_k,\theta_{k+1}]$\,,\,$k=1,\ldots,N-1$ and $[\theta_N,1]\cup[0,\theta_1]$,
\item If $e_k$ is an elliptic singular point and if $|\theta-\theta_k|$ is sufficiently small then we have $b(\theta)=-\frac{1}{2}a(\theta)(\theta-\theta_k).$
\end{itemize}
\end{proposition}
We consider first the situation near boundary singular points.

\subsection{A normal form for the spanning surface near boundary
singularities} 
We first simplify the surface ${\mathcal D}$ near singular points on the boundary:
\begin{proposition}\label{3.1.2.}
Let ${\mathcal L}$ be a Legendrian knot in a three dimensional
contact manifold $(M,\lambda)$ and let ${\mathcal D}\subset M$ be
an embedded surface with $\partial {\mathcal D}={\mathcal L}$.
Assume that the singular points of the characteristic foliation on
${\mathcal D}$ are nondegenerate. Denote by $(\theta,x,y)\in
S^1\times{\bf R}$ the coordinates near ${\mathcal L}$ provided by
the Legendrian neighborhood theorem. If ${\mathcal D}$ is
parameterized by
\[\{(\theta,x(\theta,r),y(\theta,r))\in S^1\times{\bf
R}\,|\,r\in [0,1]\}\] near ${\mathcal L}$ then there is an
embedded surface ${\mathcal D}'$ with the following properties:
\begin{itemize}
\item ${\mathcal D}'$ is obtained from ${\mathcal D}$ by a
$C^0$--small perturbation supported near the boundary singular
points leaving the boundary fixed, i.e. $\partial{\mathcal
D}'=\partial{\mathcal D}={\mathcal L}$.
\item ${\mathcal D}'$ has the same singular points as ${\mathcal
D}$.
\item If $(\theta_0,0,0)$ is a boundary singularity and
\[{\mathcal D}'=\{(\theta,x'(\theta,r),y'(\theta,r))\in
S^1\times{\bf R}\,|\,r\in[0,1]\}\] then
\[y'(\theta,r)=cx'(\theta,r)(\theta-\theta_0)+\frac{b}{2}x'(\theta,r)^2\]
where
\begin{enumerate}
\item $c=-\frac{1}{2}$ and $b=0$ if $(\theta_0,0,0)$ is elliptic,
\item $c\in (-\infty,-1)\cup (0,+\infty)$ if $(\theta_0,0,0)$ is
hyperbolic.
\end{enumerate}
\end{itemize}
\end{proposition}
{\bf Proof:}\\ Let us first point out how to recognize the type of
the singularity $(\theta_0,0,0)$ in the above parameterization.
Since the Jacobian of the map
$\Psi(\theta,r)=(\theta,x(\theta,r))$ at the point $(\theta_0,0)$
has rank $2$, there is a local inverse and we parameterize
${\mathcal D}$ by
\[{\mathcal D}=\{(\theta,x,(y\circ\Psi^{-1})(\theta,r))\},\]
where $(\theta,x)$ is sufficiently near to $(\theta_0,0)$. Note
that in the case of a positive (negative) singular point
$(\theta_0,0,0)$ the map $\Psi^{-1}$ is only defined for
non-positive (non-negative) $x$. We write
\[f(\theta,x):=(y\circ\Psi^{-1})(\theta,x)\]
and note that
\begin{itemize}
\item $f(\theta,0)\equiv 0$,
\item $\partial_xf(\theta_0,0)=0$ since $(\theta_0,0,0)$ is a
singular point.
\end{itemize}
We may extend $f$ smoothly so that it is defined for small $|x|$
regardless of the sign of $x$. Write
\[f(\theta,x)=\frac{a}{2}(\theta-\theta_0)^2+\frac{b}{2}x^2+c(\theta-\theta_0)x+h(\theta,x)\]
with $a=\partial_{\theta\theta}f(\theta_0,0)$ ,
$b=\partial_{xx}f(\theta_0,0)$ ,
$c=\partial_{x\theta}f(\theta_0,0)$, and $h$ of order at least $3$
in $(\theta-\theta_0,x)$. Note that $a=0$ and also
$h(\theta,0)=0$, hence
\[
f(\theta,x)=\frac{b}{2}x^2+c(\theta-\theta_0)x+h(\theta,x).
\]
Investigate now the admissible values for the constants $b$ and
$c$. The surface ${\mathcal D}$ is given by $H^{-1}(0)$, where
\[
H(\theta,x,y):=y-f(\theta,x).
\]
Then the vector field $\hat{V}_H$, which is defined by
$i_{\hat{V}_H}d\lambda=(i_{X_{\lambda}}dH)d\lambda-dH$ and
$i_{\hat{V}_H}\lambda=0$, is given by
\[
\hat{V}_H(\theta,x,y)=-\partial_xf(\theta,x)\frac{\partial}{\partial\theta}+(x+\partial_{\theta}
f(\theta,x))\frac{\partial}{\partial
x}+x\partial_xf(\theta,x)\frac{\partial}{\partial y},
\]
and it induces the characteristic foliation on ${\mathcal D}$.
Its linearization $\hat{V}'_H(\theta,x,y)$ is given by
\[
\left(\begin{array}{ccc} -\partial_{x\theta}f(\theta,x) &
-\partial_{xx}f(\theta,x) & 0 \\
\partial_{\theta\theta}f(\theta,x) &
1+\partial_{x\theta}f(\theta,x) & 0 \\
x\partial_{x\theta}f(\theta,x) &
\partial_xf(\theta,x)+x\partial_{xx}f(\theta,x) & 0
\end{array}\right).
\]
The contact structure $\ker\lambda(\theta_0,0,0)$ is generated by
the vectors $(1,0,0)$ and $(0,1,0)$. We represent
$\hat{V}'_H(\theta_0,0,0)$ by the matrix
\[
\left(\begin{array}{cc} -\partial_{x\theta}f(\theta_0,0) &
-\partial_{xx}f(\theta_0,0) \\ 0 &
1+\partial_{x\theta}f(\theta_0,0)
\end{array}\right)=
\left(\begin{array}{cc} -c & -b \\ 0 & 1+c \end{array}\right).
\]
The singular point $(\theta_0,0,0)$ is then hyperbolic if $c(c+1)>0$ and elliptic if $c(c+1)<0$. Let us
translate this into our original parameterization
$(\theta,x(\theta,r),y(\theta,r))$ of ${\mathcal D}$. Using
$f\circ\Psi=y$, we compute
\[
\partial_{\theta r}y(\theta_0,0)=c\,\partial_r x(\theta_0,0),
\]
so that we are in the following situation: A singular point
$(\theta_0,0,0)$ is
\begin{enumerate}
\item positive if $\partial_rx(\theta_0,0)<0$,
\item negative if $\partial_rx(\theta_0,0)>0$,
\item elliptic if $\frac{\partial_{\theta
r}y}{\partial_rx}(\theta_0,0)\in (-1,0)$ and
\item hyperbolic if $\frac{\partial_{\theta
r}y}{\partial_rx}(\theta_0,0)\in (-\infty,-1)\cup (0,+\infty).$
\end{enumerate}
We will remove now the higher order term $h$ by a perturbation.
Take a smooth function $\beta:[0,\infty)\rightarrow[0,1]$ with
$\beta\equiv 0$ on $[0,1]$, $\beta\equiv 1$ on $[2,\infty)$ and
$\beta'\ge 0$. For small $\delta>0$ we define
$\beta_{\delta}:=\beta\left(\frac{\theta^2+x^2}{\delta^2}\right)$
and
\[
f_{\delta}(\theta,x):=\frac{b}{2}x^2+cx(\theta-\theta_0)+\beta_{\delta}
(\theta-\theta_0,x)h(\theta,x).
\]
This perturbation takes place in a small neighborhood of the
singular point $(\theta_0,0,0)$. We have to show that the new
surface given by the graph of $f_{\delta}$ has the same
singularities as ${\mathcal D}$ provided $\delta>0$ was chosen
sufficiently small. We proceed indirectly. Assume that for any
sequence $\delta_n\searrow 0$ there is a singular point
$(\theta_n,x_n,f_{\delta_n}(\theta_n,x_n))$ on the surface
${\mathcal D}_{\delta_n}$ given by the graph of $f_{\delta_n}$
which satisfies $(\theta_n-\theta_0)^2+x_n^2\le2\delta_n^2$ and is
different from $(\theta_0,0,0)$. If
$(\theta_n,x_n,f_{\delta_n}(\theta_n,x_n))$ is singular then
\begin{eqnarray*}
0 & = & \partial_{\theta}f_{\delta_n}(\theta_n,x_n)+x_n \\
 & = &
 (c+1)x_n+\beta_{\delta_n}(\theta_n-\theta_0,x_n)\partial_{\theta}
 h(\theta_n,x_n)+ \\
  & &
  +\frac{2(\theta_n-\theta_0)}{\delta^2_n}\beta'\left(\frac{(\theta_n-\theta_0)^2+x_n^2}{
  \delta^2_n}\right)h(\theta_n,x_n)
\end{eqnarray*}
and
\begin{eqnarray*}
0 & = & \partial_xf_{\delta_n}(\theta_n,x_n) \\
 & = &
 bx_n+c(\theta_n-\theta_0)+\beta_{\delta_n}(\theta_n-\theta_0,x_n)\partial_xh(\theta_n,x_n)
 + \\
  & &
  +\frac{2x_n}{\delta^2_n}\beta'\left(\frac{(\theta_n-\theta_0)^2+x_n^2}{\delta_n^2}\right)h(\theta_n,
  x_n).
\end{eqnarray*}
We write shortly
\begin{eqnarray}\label{normalform1}
0 & = &
(c+1)x_n+\beta_{\delta_n}\partial_{\theta}h+\frac{2(\theta_n-\theta_0)}{\delta^2_n}\beta'_{\delta_n}
h \\ 0 & = &
bx_n+c(\theta_n-\theta_0)+\beta_{\delta_n}\partial_xh+\frac{2x_n}{\delta^2_n}\beta'_{\delta_n}h.
\nonumber
\end{eqnarray}
{\bf Remark:}\\
The reader should be aware that $\beta'_{\delta_n}$ is not the derivative of $\beta_{\delta_n}$, but the rescaled derivative of $\beta$
\[
\beta'_{\delta_n}:=\beta'\left(\frac{(\theta_n-\theta_0)^2+x_n^2}{\delta_n^2}\right).
\]
\newline
Equation (\ref{normalform1}) is the same as
\begin{equation}\label{normalform1a}
\left(\begin{array}{c} x_n \\ \theta_n-\theta_0 \end{array}\right)=\frac{1}{c(c+1)}\left(\begin{array}{cc} c & 0 \\ -b & c+1 \end{array}\right)H(\theta_n,x_n)
\end{equation}
with
\[
H(\theta_n,x_n)=-\left(\begin{array}{c} \beta_{\delta_n}\partial_{\theta}h+\frac{2(\theta_n-\theta_0)}{\delta^2_n}\beta'_{\delta_n}
h \\ \beta_{\delta_n}\partial_xh+\frac{2x_n}{\delta^2_n}\beta'_{\delta_n}h\end{array}\right),
\]
which satisfies
\[
\frac{H(\theta_n,x_n)}{\sqrt{(\theta_n-\theta_0)^2+x_n^2}}\rightarrow 0
\]
as $n\rightarrow \infty$ since $H$ is of order at least $2$ in $(\theta_n-\theta_0,x_n)$. Dividing equation (\ref{normalform1a}) by $\sqrt{(\theta_n-\theta_0)^2+x_n^2}$ and passing to the limit $n\rightarrow\infty$ we obtain a contradiction.\\ 
Hence we may assume that ${\mathcal D}$ is given
by the graph of
\[
f(\theta,x)=\frac{b}{2}x^2+cx(\theta-\theta_0)
\]
for $x^2+(\theta-\theta_0)^2$ sufficiently small. Consider now the
case where $(\theta_0,0,0)$ is an elliptic singularity. Take the
same smooth function $\beta$ as before and define for $\delta>0$
\[
f_{\delta}(\theta,x):=\beta\left(\frac{(\theta-\theta_0)^2+x^2}{\delta^2}\right)\frac{b}{2}x^2+
cx(\theta-\theta_0),
\]
so that $f_{\delta}\equiv f$ if $(\theta-\theta_0)^2+x^2\ge
2\delta^2$ and $f_{\delta}(\theta,x)=cx(\theta-\theta_0)$ if
$(\theta-\theta_0)^2+x^2\le\delta^2$. Writing
$\beta_{\delta}:=\beta\left(\frac{(\theta-\theta_0)^2+x^2}{\delta^2}\right)\,
,\,
\beta'_{\delta}:=\beta'\left(\frac{(\theta-\theta_0)^2+x^2}{\delta^2}\right)$ as before, the condition of $(\theta,x,f_{\delta}(\theta,x))$ being a singular point is
\[
0 = \left(\begin{array}{cc} c & b\beta_{\delta}
+\frac{bx^2}{\delta^2}\beta'_{\delta} \\
\frac{bx^2}{\delta^2}\beta'_{\delta} & c+1
\end{array}\right)\left(\begin{array}{c} \theta-\theta_0 \\ x
\end{array}\right)
\]
which implies
\begin{eqnarray*}
0 & = &
c(c+1)-b^2\left(\frac{x^2}{\delta^2}\beta_{\delta}\beta'_{\delta}+\frac{x^4}{\delta^4}(\beta'_{\delta})^2
\right) \\
 & \le & c(c+1)
 \end{eqnarray*}
in contradiction to the fact that $(\theta_0,0,0)$ is an elliptic
singularity. Hence we may assume that
\[
f(\theta,x)=cx(\theta-\theta_0)
\]
near an elliptic singularity $(\theta_0,0,0)$, where $c\in
(-1,0)$. Now we will carry out a last modification to achieve
$c=-\frac{1}{2}$. We take a smooth function
\[
\beta:{\bf R}\longrightarrow
[\mbox{min}\{c,-1/2\},\mbox{max}\{c,-1/2\}]\subset (-1,0)
\]
with $\beta(s)=-\frac{1}{2}$ for $|s|\le 1$ and $\beta(s)=c$ for
$|s|\ge 2$. Define for small $\delta>0$
\[
f_{\delta}(\theta,x):=\beta\left(\frac{(\theta-\theta_0)^2+x^2}{\delta^2}\right)x(\theta-\theta_0).
\]
Again, we did not create any new singular points. This completes
the proof of proposition \ref{3.1.2.}.\qed \\

\subsection{Perturbing the spanning surface near the
Legendrian knot}

We will now show the following:
\begin{proposition}\label{3.1.3.}
Let ${\mathcal L}$ be a Legendrian knot in a three dimensional
contact manifold $(M,\lambda)$ and let ${\mathcal D}$ be an
embedded surface with $\partial{\mathcal D}={\mathcal L}$. Assume
that the singular points on ${\mathcal D}$ are non--degenerate.\\ Then there is an embedded surface
${\mathcal D}'$ with $\partial{\mathcal D}'={\mathcal L}$ which
differs from ${\mathcal D}$ by a $C^0$--small perturbation
supported near ${\mathcal L}$ and leaving ${\mathcal L}$ fixed so
that ${\mathcal D}$ and ${\mathcal D}'$ have the same singular
points and the following holds:\\ There is a neighborhood $U$ of
${\mathcal L}$ and a diffeomorphism $\Phi:U\rightarrow
S^1\times{\bf R}^2$ so that
\begin{enumerate}
\item $\Phi^{\ast}(dy+xd\theta)=\lambda|_U$ , $(\theta,x,y)\in
S^1\times {\bf R}$,
\item $\Phi({\mathcal L})=S^1\times\{(0,0)\}$,
\item $\Phi(U\cap{\mathcal
D}')=\{(\theta,a(\theta)r,b(\theta)r)\in S^1\times{\bf
R}\,|\,\theta\in S^1\approx{\bf R}/{\bf Z}\, ,\, r\in[0,1]\}$,
where $\theta\longmapsto\left(\begin{array}{c} a(\theta) \\
b(\theta)
\end{array}\right)$ is a smooth closed curve in ${\bf
R}^2\backslash\{0\}$ with the following properties:
\begin{enumerate}
\item $\mbox{tb}({\mathcal L})=\,\mbox{deg}\left[\theta\longmapsto\left(\begin{array}{c}
a(\theta) \\ b(\theta)
\end{array}\right)\right]$, where $tb$ denotes denotes the Thurston--Bennequin invariant of the Legendrian knot (see \cite{Eliashberg}).
\item $b(\theta_0)=0$ if and only if $\Phi^{-1}(\theta_0,0,0)$ is
a singular point on ${\mathcal L}$,
\item a singular point $\Phi^{-1}(\theta_0,0,0)$ is
\begin{enumerate}
\item positive (negative) if $a(\theta_0)<0\ (a(\theta_0)>0)$,
\item elliptic if $c=\frac{b'(\theta_0)}{a(\theta_0)}\in (-1,0)$,
\item hyperbolic if $c=\frac{b'(\theta_0)}{a(\theta_0)}\in
(-\infty,-1)\cup(0,+\infty)$,
\end{enumerate}
\item for $\theta$ near $\theta_0$, where $\Phi^{-1}(\theta_0,0,0)$
is a singular point, we have
$b(\theta)=c\,(\theta-\theta_0)\,a(\theta)$,
\item if $\Phi^{-1}(\theta_0,0,0)$ and $\Phi^{-1}(\theta_1,0,0)$
are singular points of opposite sign with
$\theta_0<\theta_1$, so that all the points $(\theta,0,0)$ are
non--singular for $\theta\in(\theta_0,\theta_1)$, then $a$ has
exactly one zero in the interval $(\theta_0,\theta_1)$.
\end{enumerate}
\end{enumerate}
\end{proposition}
{\bf Proof:}\\ We parameterize ${\mathcal D}$ again by
\[
{\mathcal D}=\{(\theta,x(\theta,r),y(\theta,r))\in S^1\times{\bf
R}^2\,|\,\theta\in S^1\, ,\,r\in[0,1]\}
\]
and we expand $x,y$ as follows:
\begin{eqnarray*}
x(\theta,r) & = & \partial_rx(\theta,0)\,r+h(\theta,r) \\
y(\theta,r) & = & \partial_ry(\theta,0)\,r+k(\theta,r),
\end{eqnarray*}
where $h,k$ are of order at least $2$ in $r$ and $1$--periodic in
$\theta$. For small $r$ and $|\theta-\theta_0|$, where
$(\theta_0,0,0)$ is a boundary singular point, we have
\begin{equation}\label{normalform3}
y(\theta,r)=c\,x(\theta,r)(\theta-\theta_0)+\frac{b}{2}x^2(\theta,r)
\end{equation}
by proposition \ref{3.1.2.}. In the case of an elliptic
singularity we may assume that $c=-\frac{1}{2}$ and $b=0$. We want
to perturb ${\mathcal D}$ near its boundary leaving
$\partial{\mathcal D}$ fixed, so that the higher order terms $h$
and $k$ disappear. We will only indicate the necessary steps and leave the details to the reader. The verification that no new singularities are created is completely straight forward using the normal form (\ref{normalform3}) near the singular points. Pick a smooth function
$\beta:[0,\infty)\rightarrow[0,1]$ with $\beta\equiv 0$ on
$[0,1]$, $\beta\equiv 1$ on $[2,\infty)$ and $0\le\beta'(s)\le 2$
for all $s\ge 0$. We define
\begin{eqnarray*}
x_{\delta}(\theta,r) &
:= & \partial_rx(\theta,0)r+\beta\left(\frac{r}{\delta}\right)h(\theta,r),\\
y_{\delta}(\theta,r) & := &
\partial_ry(\theta,0)r+\beta\left(\frac{r}{\delta}\right)k(\theta,r)
\end{eqnarray*}
and
\[
{\mathcal
D}_{\delta}=\{(\theta,x_{\delta}(\theta,r),y_{\delta}(\theta,r))\in
S^1\times{\bf R}^2\,|\,\theta\in S^1\,,\,r\in[0,1]\}.
\]
For $r\ge 2\delta$ the perturbed surface ${\mathcal D}_{\delta}$
coincides with ${\mathcal D}$ and we have $\partial{\mathcal
D}_{\delta}=\partial{\mathcal D}={\mathcal L}$. The surface ${\mathcal D}_{\delta}$ has the same singularities on the
boundary as ${\mathcal D}$ and that ${\mathcal D}_{\delta}$ has no
singularities in the range $0<r<2\delta$ provided $\delta>0$ was
chosen sufficiently small. It remains to verify that the surface ${\mathcal D}_{\delta}$ is embedded for sufficiently small $\delta$. If it were not then we could find sequences $\delta_k\searrow 0$ , $0\le r_k,r'_k\le 2\delta_k$ and $\theta_k$ such that 
\[
(\partial_rx_{\delta_k}(\theta_k,r_k),\partial_r y_{\delta_k}(\theta_k,r_k))=(0,0)
\]
for all $k$ (surface not immersed) or
\[
x_{\delta_k}(\theta_k,r_k)=x_{\delta_k}(\theta_k,r'_k)
\]
\[
y_{\delta_k}(\theta_k,r_k)=y_{\delta_k}(\theta_k,r'_k)
\]
for all $k$ (surface has self--intersections). Both assertions contradict the fact that $(\partial_rx(\theta,0),\partial_ry(\theta,0))\neq (0,0)$ for all $\theta$, and can therefore not occur.\\
Hence we may assume that near ${\mathcal L}$ we have
\[
{\mathcal D}=\{(\theta,a(\theta)r,b(\theta)r)\in S^1\times{\bf
R}^2\,|\,r\in[0,1]\},
\]
where the map $\theta\longmapsto\left(\begin{array}{c} a(\theta)
\\ b(\theta) \end{array} \right)$ is a closed curve in ${\bf
R}^2\backslash\{0\}$. A point $(\theta_0,0,0)$ is a singular point
if and only if $b(\theta_0)=0$ and it is
\begin{enumerate}
\item positive if $a(\theta_0)<0$,
\item negative if $a(\theta_0)>0$,
\item elliptic if $\frac{b'(\theta_0)}{a(\theta_0)}\in (-1,0)$ and
\item hyperbolic if $\frac{b'(\theta_0)}{a(\theta_0)}\in
(-\infty,-1)\cup(0,+\infty)$.
\end{enumerate}
If $r$ and $\theta-\theta_0$ are sufficiently small, where
$(\theta_0,0,0)$ is a singular point, then we compute with the
normal form (\ref{normalform3}):
\begin{eqnarray}\label{normalform8}
b(\theta) & = & \partial_ry(\theta,0) \nonumber \\
 & = & c\,\partial_rx(\theta,0)(\theta-\theta_0) \\
 & = & c\,a(\theta)(\theta-\theta_0),\nonumber
\end{eqnarray}
so if we use the parameters $(\theta,\rho=a(\theta)r)$ instead of
$(\theta,r)$ then ${\mathcal D}$ is given by
\[
{\mathcal D}=\{(\theta,\rho,c\rho(\theta-\theta_0))\}\ \mbox{ near
} (\theta_0,0,0).
\]
Hence the modification that we carried out on ${\mathcal D}$ in
this section did not affect the normal form near boundary
singularities that we have constructed in the previous section.\\
In this picture it is easy to understand the Thurston--Bennequin
invariant of the knot ${\mathcal L}$. Let us shift ${\mathcal L}$
along the Reeb vector field to get a knot
\[
{\mathcal L}':=\{(\theta,0,\delta)\in S^1\times {\bf R}^2
\,|\,\theta\in[0,1]\}
\]
with some small constant $\delta$. Then ${\mathcal L}'$ and
${\mathcal D}$ intersect if and only if
\[
a(\theta)=0
\]
and
\[
r=\frac{\delta}{b(\theta)}.
\]
The condition $a(\theta)=0$ means that the Reeb vector field
$X_{\lambda}$ is tangent to ${\mathcal D}$ at the point
$(\theta,0,\delta)$. Without affecting the value of the
intersection number  int$({\mathcal L}',{\mathcal D})$ we may
perturb the loop $(a(\theta),b(\theta))$ slightly so that
$a'(\theta)\neq 0$ whenever $a(\theta)=0$. Then we compute with
$\lambda=dy+xd\theta$ and
\[
{\mathcal S}:=\{\theta\in[0,1]\,|\,a(\theta)=0\mbox{ and
sign}\,b(\theta)=\mbox{ sign}\,\delta\}:
\]
\begin{eqnarray}\label{normalform9}
\mbox{tb}({\mathcal L}) & = & \sum_{\theta\in{\mathcal
 S}}\mbox{sign}\left[(\lambda\wedge
 d\lambda)_{(\theta,0,\delta)}\left(
 \left(\begin{array}{c} 1 \\ \frac{a'(\theta)\delta}{b(\theta)} \\
 \frac{b'(\theta)\delta}{b(\theta)} \end{array}\right),
 \left(\begin{array}{c} 0 \\ 0 \\ b(\theta)  \end{array}\right),
 \left(\begin{array}{c} 1 \\ 0 \\ 0  \end{array}\right)
\right)\right] \nonumber \\
 & = & \sum_{\theta\in{\mathcal S}}\mbox{sign}[-\delta a'(\theta)]
 \\
 & = & \mbox{deg}\left[\theta\longmapsto\left(\begin{array}{c}
 a(\theta) \\ b(\theta) \end{array}\right)\right].\nonumber
\end{eqnarray}
Assume now that $(\theta_0,0,0)$ and $(\theta_1,0,0)$ are singularities of opposite sign with $\theta_0<\theta_1$, so that
all the points $(\theta,0,0)$ with $\theta\in(\theta_0,\theta_1)$
are not singular. Let us assume that $(\theta_0,0,0)$ is the
negative singularity. Then
\begin{itemize}
\item $b(\theta_0)=b(\theta_1)=0$ and $b$ is nonzero on
$(\theta,\theta_1)$.
\item $a(\theta_0)>0\ ,\ a(\theta_1)<0.$
\end{itemize}
We would like to perturb ${\mathcal D}$ near ${\mathcal L}$,
leaving the boundary fixed, so that $a$ has only one zero in the
interval $(\theta_0,\theta_1)$. Let $\delta>0$ and pick a smooth
function $\beta$ so that $\beta\equiv 0$ on $[0,\delta]$ and
$\beta\equiv 1$ on $[2\delta,\infty)$. Let $\hat{a}$ be a
1--periodic function which coincides with $a$ except on some interval $[\theta_0+\ve,\theta_1-\ve]$, and which has exactly one zero between $\theta_0$ and
$\theta_1$. We define
\[
\tilde{a}(\theta,r):=(1-\beta(r))\hat{a}(\theta)+\beta(r)a(\theta)
\]
and denote the new surface by
\[
{\mathcal D}_{\delta}:=\{(\theta,\tilde{a}(\theta,r)r,b(\theta)r)\},
\]
which has the same number and type of singularities as ${\mathcal
D}$ because we did not change the function $b$ and because
$\tilde{a}$ coincides with $a$ near $\theta_0$ and $\theta_1$.
Moreover, ${\mathcal D}_{\delta}$ is embedded since it is immersed
and the map
$(\theta,r)\longmapsto(\theta,\tilde{a}(\theta,r),b(\theta))$ is
also injective. This completes the proof of proposition
\ref{3.1.3.}.\qed \\

{\bf Remarks:}
The negative singularities correspond to the points where
the curve
$\theta\stackrel{\gamma}{\longmapsto}a(\theta)+ib(\theta)\in{\bf
C}\backslash\{0\}$ hits the positive real axis. Similarly,
positive singularities corresponds to the intersection of $\gamma$
with the negative real axis.

\subsection{The non--Lagrangian part of the boundary
condition}\label{3.2.}
The submanifold $\Real\times{\mathcal L}$ is a Lagrangian
 submanifold of the symplectisation $(\RM, d(e^t\lambda))$.
 However, the submanifold $\{0\}\times{\mathcal D}$ is
 only totally real with respect to any $\tilde{J}$ away from the singular points. These two submanifold serve as boundary conditions for our boundary value problem, and
 we have to find a way to deal with $\{0\}\times{\mathcal D}$ in
 order to derive apriori estimates. The problem is the following: The fibers of the vector bundles $T({\bf R}\times{\mathcal L})$ and $\tilde{J}T({\bf R}\times{\mathcal L})$ are orthogonal with respect to the $\tilde{J}$--invariant metric $g=d(e^t\lambda)\circ(\tilde{J}\times\mbox{ Id})$ while $T(\{0\}\times{\mathcal D})$ and $\tilde{J}T(\{0\}\times{\mathcal D})$ are only transverse, but not orthogonal. On the other hand, we will need this orthogonality to prove asymptotic decay estimates later (without it certain operators would fail to be self--adjoint). The way out is the following: Instead of using the metric $g$ above, we use a different one where we have orthogonality. We will be able to control this metric if we do estimates later on.
There is a 2--form $\omega$ near the intersection set $\{0\}\times {\mathcal L}$ of ${\bf R}\times {\mathcal L}$ and $\{0\}\times {\mathcal D}$ which is nondegenerate away from the singular points so that both submanifolds become Lagrangian with respect to $\omega$, and $\omega$ is compatible with the almost complex structure $\tilde{J}$. In general, we cannot expect $\omega$ to be closed, unless we weaken our requirements and replace compatibility by tameness (i.e. $\omega(v,\tilde{J}v)>0$ for all $v\neq 0$). It will turn out that we need the compatibility condition, but we do not need $\omega$ to be closed. We construct such a 2--form explicitly in local coordinates. We will confine ourselves to a special almost complex structure $\tilde{J}$ near $\{0\}\times{\mathcal L}$ which will also be used in the subsequent papers \cite{part2} , \cite{part3} and \cite{part4}.\\

{\sl From now on we pick an almost complex structure $\tilde{J}$ on $\RM$, where the corresponding $J:\xi\rightarrow\xi$ has the following form in local coordinates} near $\{0\}\times{\mathcal L}$: 
\begin{equation}\label{local-formula-for-J}
J(\theta,x,y)\cdot(1, 0, -x):= (0,-1,0)\ ,\ J(\theta,x,y)\cdot(0,1,0):=(1,0,-x).
\end{equation}

\begin{lemma}\label{deal_with_totally_real}
 If $U_k\subset M$ are disjoint open neighborhoods of the singular points $e_k$ , $k=1,\ldots,N$ on the boundary ${\mathcal L}=\partial{\mathcal D}$
 then there exist an open neighborhood $V\subset M$ of ${\mathcal L}$ and
 a nondegenerate 2--form $\omega$ defined
 on $W={\bf R}\times (V\backslash \cup_k U_k)\subset \RM$,
 so that $\omega|_{T(\{0\}\times{\mathcal
 D})}\equiv 0$ , $\omega|_{T(\Real\times{\mathcal L})}\equiv 0$
 and the form $\omega$ is compatible with $\tilde{J}$, i.e.
 $\omega\circ(\mbox{Id}\times\tilde{J})$ is a Riemannian metric.
 \end{lemma}

{\bf Proof of lemma \ref{deal_with_totally_real}:}\\
 Use the coordinates $(\theta,x,y)\in\Real^3$ near
${\mathcal L}$
 which we derived in section \ref{normal-form-for-the-surface}, where the contact
 form equals $dy+xd\theta$ and $\{0\}\times{\mathcal D}$ is represented by
\[ 
\{(0,\theta,a(\theta)r,b(\theta)r)\in\{0\}\times{\Real}^3\,|\,r,\theta\in[0,1]\}
\]
Denoting the standard Euclidean product on ${\bf R}^4$ by $\langle\,.\,,\,.\,\rangle$, we have to find a function with values in the set of skew--symmetric $4\times 4$--matrices $\Omega(\tau,\theta,x,y)$ such that
\begin{enumerate}
\item $\langle\,.\,,\,\Omega\tilde{J}\,.\,\rangle$ is a metric,
\item $\langle v,\Omega w\rangle|_{(\tau,\theta,0,0)}=0$ for all $v,w\in T_{(\tau,\theta,0,0)}({\bf R}\times{\mathcal L})$,
\item $\langle v,\Omega w\rangle|_{(0,\theta,q(\theta)y,y)}=0$ for all $v,w\in T_{(0,\theta,q(\theta)y,y)}(\{0\}\times{\mathcal D})$, where $q(\theta):=a(\theta)/b(\theta)$.
\end{enumerate}
The matrix of $\tilde{J}$ is given by
\[
\tilde{J}(\tau,\theta,x,y)=\left(\begin{array}{cccc} 0 & -x & 0 & -1 \\ 0 & 0 & 1 & 0 \\ 0 & -1 & 0 & 0 \\ 1 & 0 & -x & 0 \end{array}\right).
\]
We write $\Omega=(\omega_{kl})_{1\le k,l\le 4}$ with $\omega_{kl}=-\omega_{lk}$. If we choose
\begin{eqnarray*}
\omega(\tau,\theta,x,y) & = & -xC\,d\tau\wedge d\theta-q(\theta)d\tau\wedge dx+\\
 & & +C d\tau\wedge dy-d\theta\wedge dx+q(\theta)d\theta\wedge dy
\end{eqnarray*}
and $\omega=\langle\,.\,,\Omega\,.\,\rangle\,$, where
\[
C>\mbox{max}\{0,\sup q^2(\theta)\},
\]
then the matrix $\Omega\tilde{J}$ is positive definite if $x,y$ are sufficiently small. 
\qed \\

\section{Local existence of solutions}\label{local-existence-of-sol}

In this section we establish local fillings by pseudoholomorphic curves near an elliptic singularity at the boundary.
 Because of proposition \ref{3.1.1.} we are in the following situation near an elliptic singular point $e\in{\mathcal L}$:\\
We may assume that the contact manifold is the three dimensional
Euclidean space $\{(\theta,x,y)\in{\bf R}^3\}$ endowed with the
contact form $\lambda=dy+xd\theta$. The piece of the Legendrian
knot situated near $e$ corresponds to some interval
$\{(\theta,0,0)\in{\bf R}^3\,|\,|\theta|<\varepsilon\}$, where
$\varepsilon>0$ is a suitable constant. The elliptic
singular point then corresponds to the origin in ${\bf R}^3$ and the spanning surface
${\mathcal D}$ is given by $\{(\theta,x,-\frac{1}{2}\theta
x)\in{\bf R}^3\,|\,|\theta|<\varepsilon\,,\,x\le 0\}$ if $e$ is a positive elliptic point, otherwise we have $x\ge 0$. We start constructing solutions near $e$. The contact structure is generated by the
vectors
\[e_1=\left(\begin{array}{c} 1 \\ 0 \\ -x \end{array}\right)\ \mbox{and}\ e_2=\left(\begin{array}{c} 0 \\ 1 \\ 0 \end{array}\right).\]
Recall that we have chosen a particular complex structure $J$ on $\ker\lambda$ near the Legendrian knot by demanding
\begin{equation}\label{standard-J-near-e+-}
Je_1:=-e_2\ \mbox{and}\ Je_2=e_1.
\end{equation}
This complex structure is compatible with $d\lambda$, i.e.
$d\lambda\circ(\mbox{Id}\times J)$ is a bundle metric and defines an almost
complex structure $\tilde{J}$ in the usual way.\\ The
boundary value problem, we are going to study, is the following:
\[
 \begin{array}{ll}
          \tilde{u}=(a,u):S\longrightarrow\RM & \\
          \partial_s\tu+\tilde{J}(\tu)\partial_t\tu=0 & \\
          \tu(s,0)\subset{\bf R}\times{\cal L} & \\
          \tu(s,1)\subset\{0\}\times{\cal D}^{\ast} & \\
           0< E(\tu)<+\infty & 
          \end{array},
\]
where ${\mathcal D}^{\ast}$ is
the spanning surface without the singular points, $S:={\bf
R}\times[0,1]$ and
\[E(\tu)=\sup_{\phi\in\Sigma}\int_S\tu^{\ast}d(\phi\lambda)\
\mbox{ (energy of $\tu$ ) }\] with $\Sigma:=\{\phi\in
C^{\infty}({\bf R},[0,1])\,|\,\phi'\ge 0\}$. Since we have chosen
good coordinates near the elliptic singular points and an explicit
almost complex structure $\tilde{J}$, we will be able to explicitly state 
1--parameter families of solutions to the above boundary value
problem near the elliptic singular points.\\ These solutions look simpler after having performed a biholomorphic transformation of the domain as
follows: Let $\Omega:=\{z=s+it\in{\bf C}\,|\,s^2+t^2\le 1\,,\,t\ge
0\}\backslash\{-1,+1\}$ be the upper half disk in the complex
plane without the corner points. The infinite strip $S$ and
$\Omega$ are equivalent via the biholomorphic map
\[S\longrightarrow\Omega\]
\begin{equation}\label{strip-halfdisk}
s+it\longmapsto\frac{e^{\frac{\pi}{2}(s+it)}-1}{e^{\frac{\pi}{2}(s+it)}+1}=\mbox{tanh}\left(\frac{\pi}{4}(s+it)\right).
\end{equation}
Under this transformation, ${\bf R}\times\{0\}$ is mapped onto
$(-1,+1)$ and ${\bf R}\times\{1\}$ is mapped onto
$\{s+i\sqrt{1-s^2}\in{\bf C}\,|\,s\in (-1,+1)\}$. We write in
coordinates
\[\tilde{u}=(a,\theta,x,y):\Omega\rightarrow{\bf R}\times{\bf R}^3\]
and obtain the following boundary value problem:
\begin{eqnarray*}
\pas a-\pat y-x\pat \theta & = & 0 \\
\pas\theta+\pat x & = & 0 \\
\pas x-\pat\theta & = & 0 \\
\pat a+\pas y+x\pas \theta & = & 0 \\
x(s,0)\equiv y(s,0) & \equiv & 0 \\
a(s,\sqrt{1-s^2}) & \equiv & 0 \\
y(s,\sqrt{1-s^2}) & = & -\frac{1}{2}(x\theta)(s,\sqrt{1-s^2}).
\end{eqnarray*}
The following maps satisfy the above boundary
value problem as long as they stay in the coordinate patch near
the elliptic singular point:
\begin{equation}\label{local_solutions}
\tu_{\varepsilon}(s,t)=\left(\frac{1}{4}\varepsilon^2(s^2+t^2-1),\varepsilon
s,-\varepsilon t,\frac{1}{2}\varepsilon^2st\right),
\end{equation}
with $\ve>0$ if $e$ is a positive elliptic singular point and $\ve<0$ otherwise. 
Transforming back the infinite strip $S={\bf R}\times[0,1]$ the solutions (\ref{local_solutions}) become
\begin{eqnarray}\label{local-solutions-2}
\tu_{\ve}(s,t) & = & \left(
-\frac{\ve^2\cos\left(\frac{\pi t}{2}\right)}{2\left[\cos\left(\frac{\pi t}{2}\right)+\cosh\left(\frac{\pi s}{2}\right)\right]}, \frac{\ve \sinh\left(\frac{\pi s}{2}\right)}{\cos\left(\frac{\pi t}{2}\right)+\cosh\left(\frac{\pi s}{2}\right)}, \right. \nonumber \\
 & & \left.\frac{-\ve\sin\left(\frac{\pi t}{2}\right)}{\cos\left(\frac{\pi t}{2}\right)+\cosh\left(\frac{\pi s}{2}\right)}, \frac{\ve^2\sin\left(\frac{\pi t}{2}\right)\sinh\left(\frac{\pi s}{2}\right)}{2\left[\cos\left(\frac{\pi t}{2}\right)+\cosh\left(\frac{\pi s}{2}\right)\right]^2}
\right).
\end{eqnarray}

\section{Asymptotic behavior at
infinity}\label{asymptotic-behaviour} Assume we have a solution
of:
\begin{equation}\label{5.1.1.}
\left\{\begin{array}{ll}
          \tilde{u}=(a,u):S\longrightarrow\RM & \\
          \partial_s\tu+\tilde{J}(\tu)\partial_t\tu=0 & \\
          \tu(s,0)\subset{\bf R}\times{\cal L} & \\
          \tu(s,1)\subset\{0\}\times{\cal D}^{\ast\ast} & \\
           E(\tu)<+\infty
          \end{array}\right.
\end{equation}
where $S:={\bf R}\times[0,1]$ and ${\cal D}^{\ast\ast}$ is the
spanning surface ${\cal D}$ without some open neighborhood $U$
of the set of singular points $\Gamma$. We will show that the
condition of finite energy forces the solution to converge to
points on the knot ${\cal L}$ for $|s|\rightarrow\infty$, more
precisely $$ \tu(s,t)\longrightarrow\tilde{p}_{\pm}
\in\{0\}\times{\cal L} $$ as $s\rightarrow\pm\infty$ uniformly in
$t$. We will also show that this convergence is of exponential
nature. This fact will be crucial for the nonlinear Fredholm theory in \cite{part2}.

\subsection{The solutions approach the Legendrian asymptotically}

As a first step, we will show that the ends of a finite
energy strip $\tilde{u}$ have to approach the knot $\{0\}\times
{\mathcal L}\subset{\bf R}\times M$ asymptotically. This actually
works under the weaker assumption $\tu(s,1)\in\{0\}\times{\mathcal
D}$. The main result of this section is proposition \ref{prop5.4.} below.
\begin{lemma}\label{lemma1.1.}
Assume $\tu:S\rightarrow{\bf R}\times M$ satisfies equation
(\ref{5.1.1.}) above. If in addition $$ \int_Su^{\ast}d\lambda=0
$$ then $\tu$ must be constant.
\end{lemma}
{\bf Proof:}\\ The map $\tu=(a,u)$ satisfies the following system
of equations:
\begin{eqnarray*}
\pil\pas u+J(u)\pil\pat u & = & 0 \\ \pas a-\lambda(u)\pat u & = &
0 \\ \pat a +\lambda(u)\pas u & = & 0.
\end{eqnarray*}
Since
\begin{eqnarray*}
\int_Su^{\ast}d\lambda & = & \int_Sd\lambda(\pil\pas u,\pil\pat
u)ds\wedge dt \\
 & = & \frac{1}{2}\int_S[\,|\pil\pas u|^2_J+|\pil\pat u|^2_J\,]ds\wedge dt\\
 & = & 0,
\end{eqnarray*}
where $|.|^2_J=d\lambda(.,J.)$, we conclude that $\pil\pas
u=\pil\pat u\equiv 0$ and therefore
\begin{eqnarray*}
\triangle a\, ds\wedge dt & = & -d(da\circ i) \\
 & = & u^{\ast}d\lambda \\
 & = & 0,
\end{eqnarray*}
hence $a:S\rightarrow{\bf R}$ is harmonic and satisfies
$a(s,1)\equiv 0$. Because of $u(s,0)\in{\cal L}$ we also have $$
\pat a(s,0)=-\lambda(u(s,0))\pas u(s,0)\equiv 0. $$ Define now
$f:S\rightarrow {\bf R}$ by $$ f(s,t):=\int_0^t\pas a(s,\tau)d\tau
$$ so that $\pat f=\pas a$ and
\begin{eqnarray*}
\pas f(s,t) & = & -\int_0^t\partial_{tt}a(s,\tau)d\tau \\
 & = & -\pat a(s,t)+\pat a(s,0) \\
 & = & -\pat a(s,t).
\end{eqnarray*}
Then $\Phi:=a+if:S\rightarrow{\bf C}$ is holomorphic and satisfies
$$ \Phi(s,0)\in{\bf R} $$ $$ \Phi(s,1)\in i{\bf R}. $$ {\bf CASE
1: $|\nabla\Phi|$ is bounded}\\ We define $$
\hat{\Phi}:\hat{S}:={\bf R}\times[-1,+1]\longrightarrow{\bf C} $$
by $$ \hat{\Phi}(s,t):=\left\{\begin{array}{cc}
                        \Phi(s,t) & \mbox{if}\ t\ge 0\\
                        \overline{\Phi(s,-t)} & \mbox{if}\ t<0
                       \end{array}\right..
$$ Note that $\hat{\Phi}$ is holomorphic. Let $$
\hat{b}:=\pas(\mbox{Re}\,\hat{\Phi}):\hat{S}\longrightarrow{\bf
R}. $$ Then $\hat{b}$ is harmonic,
$C:=\sup_{\hat{S}}|\hat{b}|<+\infty$ by assumption and
$\hat{b}(s,\pm 1)\equiv 0$. Defining $$
\hat{c}(s,t):=\int^t_0\pas\hat{b}(s,\tau)d\tau-\int_0^s\pat\hat{b}(\sigma,0)d\sigma,
$$ we compute $\pat\hat{c}=\pas\hat{b}$ and
$\pas\hat{c}=-\pat\hat{b}$, hence
$\delta:=\hat{b}+i\hat{c}:\hat{S}\rightarrow{\bf C}$ is
holomorphic with bounded real part. The function $g:=e^{\delta}$
is also holomorphic and satisfies $$ |g|\le e^C\ ,\ |g(s,\pm
1)|=1. $$ Let $\eps>0$ and define a holomorphic function on
$\hat{S}$ by
\[
h_{\eps}(z):=\frac{1}{1-i\eps(z+i)}.
\]
We compute with $z=s+it$
\[
|h_{\eps}(z)|^2=\frac{1}{(1+\eps(1+t))^2+\eps^2s^2}\le 1.
\]
For $s\neq 0$ we have $|h_{\eps}(z)|^2\le\frac{1}{\eps^2s^2}$,
hence the holomorphic function $gh_{\eps}$ satisfies
\[
|g(z)h_{\eps}(z)|\le 1
\]
whenever $z\in\partial \Omega$, where
\[
\Omega:=[-\eps^{-1}e^C,\eps^{-1}e^C]\times[-1,+1].
\]
Using the maximum principle, we conclude that $|gh_{\eps}|\le 1$
on all of $\Omega$, but outside $\Omega$ we also have
\[
|g(z)h_{\eps}(z)|\le\frac{e^C}{\eps |s|}\le 1.
\]
Keeping $z\in\hat{S}$ fixed and passing to the limit $\eps\searrow
0$ we conclude that $|g(z)|=e^{\hat{b}(z)}\le 1$, hence
$\hat{b}(z)\le 0$. Repeating the same argument with $-\delta$
instead of $\delta$, we also obtain $-\hat{b}(z)\le 0$, hence
$\pas(\mbox{Re}\,\hat{\Phi})=\hat{b}(z)\equiv 0$. We know now that
$\mbox{Re}\,\hat{\Phi}$ is harmonic, does not depend on $s$ and
satisfies $\mbox{Re}\,\hat{\Phi}(s,\pm 1)\equiv 0$. This implies
that $\mbox{Re}\,\hat{\Phi}$ is identically zero and therefore
also $a\equiv 0$. In view of
\[
\pas u=\pil\pas u+(\lambda(u)\pas u)X_{\lambda}(u)
\]
and
\[
\pat u=\pil\pat u+(\lambda(u)\pat u)X_{\lambda}(u)
\]
we conclude that $u$ must be constant.\\ {\bf CASE 2:
$|\nabla\Phi|$ is unbounded}\\ Pick sequences $z'_k\in S$ ,
$\eps_k'\searrow 0$ so that
\[
\eps_k'|\nabla\Phi(z'_k)|\rightarrow +\infty.
\]
By a lemma of H. Hofer (see \cite{HZbook}, chapter 6.4, lemma 5 and \cite{AH}) we find
sequences $z_k=s_k+it_k\in S$ , $\eps_k\searrow 0$ so that
\begin{itemize}
\item $\eps_kR_k:=\eps_k|\nabla\Phi(z_k)|\longrightarrow +\infty$,
\item $|z_k-z'_k|\le\eps_k'$,
\item $|\nabla\Phi(z)|\le 2R_k$ whenever $|z-z_k|\le\eps_k$.
\end{itemize}
We may assume without loss of generality that $t_k\rightarrow
t_0\in[0,1]$. We consider the following cases after choosing a
suitable subsequence:
\begin{enumerate}
\item {$-t_kR_k\longrightarrow -\infty$
       \begin{enumerate}
         \item $R_k(1-t_k)\rightarrow l\in[0,+\infty)$
         \item $R_k(1-t_k)\rightarrow +\infty$
       \end{enumerate}}
\item $-t_kR_k\longrightarrow -l\in(-\infty,0]$, then $R_k(1-t_k)\rightarrow +\infty$
\end{enumerate}
Let us begin with the case 1b. We define
\[
\Omega_k:={\bf R}\times[-t_kR_k,R_k(1-t_k)]
\]
and the holomorphic maps $\Phi_k:\Omega_k\rightarrow{\bf C}$ by
\[
\Phi_k(z):=\Phi(z_k+zR_k^{-1})-\Phi(z_k)
\]
so that
\[
|\nabla\Phi_k(0)|=1,
\]
\[
\Phi_k(0)=0
\]
and
\[
|\nabla\Phi_k(z)|\le 2
\]
if $z\in B_{\eps_kR_k}(0)\cap\Omega_k$. Using the Cauchy integral
formula for higher derivatives we find for each compact subset $K$
of ${\bf C}$ a number $k_0$ so that $K\subset
B_{\eps_kR_k}(0)\cap\Omega_k$ for all $k\ge k_0$ and all the maps
$\Phi_k$ are bounded in $C^{\infty}(K)$ uniformly in $k\ge k_0$.
By the Ascoli--Arzela theorem, some subsequence of $(\Phi_k)$
converges in $C^{\infty}_{loc}$ to an entire holomorphic function
$\Psi$ satisfying
\[|\nabla\Psi(z)|\le 2\ ,\ |\nabla\Psi(0)|=1\ \mbox{and}\ \Psi(0)=0.\]
By Liouville's theorem $\Psi$ must be an affine function. Let
$\phi\in\Sigma$ and define $\phi_k\in\Sigma$ by
\[
\phi_k(s):=\phi(s-\mbox{Re}\,\Phi(z_k))
\]
and
\[
\tau_{\phi}(s,t):=\phi'(s)ds\wedge dt.
\]
We estimate using $u^{\ast}d\lambda=0$:
\begin{eqnarray*}
\int_{\Omega_k}\Phi_k^{\ast}\tau_{\phi} & = &
\int_S\Phi^{\ast}\tau_{\phi_k} \\
 & = & \int_S\phi'_k(a)\Phi^{\ast}(ds\wedge dt) \\
 & = & \int_S\phi'_k(a)da\wedge u^{\ast}\lambda \\
 & = & \int_S\tu^{\ast}d(\phi_k\lambda) \\
 & \le & E(\tu).
\end{eqnarray*}
For every compact $K\subset{\bf C}$ we have
\[
\int_K\Phi^{\ast}_k\tau_{\phi}\stackrel{k\rightarrow\infty}{\longrightarrow}
\int_K\Psi^{\ast}\tau_{\phi}.
\]
It follows for non constant $\phi\in\Sigma$:
\begin{eqnarray*}
+\infty & = & \int_{\bf C}\phi'(s)ds\wedge dt \\
 & = & \int_{\bf C}\tau_{\phi} \\
 & = & \int_{\bf C}\Psi^{\ast}\tau_{\phi} \\
 & \le & E(\tu).
\end{eqnarray*}
This contradiction to $E(\tu)<+\infty$ shows that case 1b. can not
occur. We will proceed similarly with the remaining cases 1a. and
2. Let us continue with case 2. We define
\[
\Omega_k:={\bf R}\times[0,R_k]
\]
and for $z=s+it\in\Omega_k$ , $z_k=s_k+it_k$
\[
\Phi_k(z):=\Phi(s_k+zR_k^{-1})-\Phi_k(s_k),
\]
so that
\[
|\nabla\Phi_k(iR_kt_k)|=1\ ,\ \Phi_k(0)=0
\]
and
\[
|\nabla\Phi_k(z)|\le 2
\]
whenever $z\in B_{\eps_kR_k}(iR_kt_k)\cap\Omega_k$. Reasoning as
before we obtain $C^{\infty}_{loc}$--convergence of some
subsequence of $(\Phi_k)$ to a holomorphic map
$\Psi:H^+\rightarrow{\bf C}$, where $H^+$ denotes the upper half
plane in ${\bf C}$. Since we have $\Phi_k({\bf R})\subset{\bf R}$
for all $k$, we also obtain
\[
\Psi(\partial H^+)\subset {\bf R}.
\]
Moreover $|\nabla\psi(z)|\le 2$ , $\Psi(0)=0$ and $\Psi$ is not
constant. Using the Schwarz reflection principle we can extend
$\Psi$ to an entire holomorphic function with bounded derivative,
so that $\Psi$ must be an affine function by Liouville's theorem.
In view of $\Psi(0)=0$ and the real boundary values we have
actually $\Psi(z)=\alpha z$ with some nonzero real number
$\alpha$. We compute as before with non constant $\phi\in\Sigma$:
\[
\int_{\Omega_k}\Phi_k^{\ast}\tau_{\phi}=\int_S\Phi^{\ast}\tau_{\phi_k}\le
E(\tu),\] where $\phi_k(s):=\phi(s-\mbox{Re}\,\Phi(s_k))$, which
implies
\[
\int_{H^+}\Psi^{\ast}\tau_{\phi}\le E(\tu)<+\infty.
\]
But on the other hand
\[
\int_{H^+}\Psi^{\ast}\tau_{\phi}=|\alpha|\int_{H^+}\phi'(s)ds\wedge
dt=+\infty,
\]
so that case 2. is impossible. We are left with case 1a. We define
$\Omega_k:={\bf R}\times[0,R_k]$ and for $z\in\Omega_k$
\[
\Phi_k(z):=\Phi(s_k+i-R_k^{-1}z)-\Phi(s_k+i).
\]
We have
\[
|\nabla\Phi_k(iR_k(1-t_k))|=1\ , \ \Phi_k(0)=0
\]
and
\[
|\nabla\Phi_k(z)|\le 2
\]
whenever $z\in B_{\eps_kR_k}(iR_k(1-t_k))\cap\Omega_k$. Moreover
\[
\Phi_k({\bf R})\subset i{\bf R}.
\]
Again, a subsequence of $(\Phi_k)$ converges in $C^{\infty}_{loc}$
to a holomorphic map
\[
\Psi:H^+\longrightarrow{\bf C}
\]
with $|\nabla\Psi(z)|\le 2$ , $\Psi(0)=0$ and $\Psi$ is not
constant. Defining
\[
\tilde{\Psi}(z):=\left\{\begin{array}{cc}
                         \Psi(z) & \mbox{if  Im}(z)\ge 0 \\
                         -\overline{\Psi(\overline{z})} & \mbox{if  Im}(z)<0
                        \end{array}\right.
\]
we obtain an entire holomorphic function with bounded gradient
which has to be affine. Since $\Psi(\partial H^+)\subset i{\bf R}$
we have $\Psi(z)=i\alpha z$ with some nonzero real number
$\alpha$. Then
\begin{eqnarray*}
E(\tu) & \ge & \int_{H^+}\tu^{\ast}d(\phi\lambda) \\
 & = & \int_{H^+}\phi'(a)da\wedge df \\
 & = & \int_{H^+}\alpha^2\phi'(-\alpha t)ds\wedge dt \\
 & = & |\alpha|\cdot\int_{H^+}\phi'(t)ds\wedge dt,
\end{eqnarray*}
but if we take a $\phi\in\Sigma$ which is not constant on
$[0,+\infty)$, then $\int_{H^+}\phi'(t)ds\wedge dt=+\infty$. This
is a contradiction to the finite energy condition. Hence we have
shown that $|\nabla\Phi|$ must be bounded, and therefore $\tu$ is
constant. \qed \\

{\bf Remarks:} There are similar results for $\tu$ defined on the
whole plane ${\bf C}$ (\cite {Hofer-Weinstein-conj}, \cite{AH}) and for $\tu$
defined on $H^{+}$ with boundary condition $\RL$ (\cite{A1}).
In the case of a finite energy strip $\tu:S\rightarrow\RM$ with
boundary condition $\tu(\partial S)\subset\RL$ we cannot conclude
from $\int_S u^{\ast}d\lambda=0$ that $\tu$ is constant (see
\cite{A1}).\\ We will omit the proof of the following lemma since it is similar to the proof of lemma \ref{lemma1.1.}:
\begin{lemma}\label{lemma1.2.}
Let $\tu=(a,u):H^+\rightarrow\RM$ be a solution of
$\pas\tu+\tilde{J}(\tu)\pat\tu=0$ satisfying the boundary
condition $\tu(\partial H^+)\subset \od$ and the finite energy
condition $E(\tu)<+\infty$. If also $\int_{H^+}u^{\ast}d\lambda=0$
then $\tu$ must be constant.
\end{lemma}
\qed
\begin{lemma}\label{lemma1.3.}
Let $\tu$ be as in equation (\ref{5.1.1.}) and assume that $u(S)$
is contained in a compact subset of $M$. Then
\[
\sup_{z\in S}|\nabla\tu(z)|<+\infty.
\]
\end{lemma}
{\bf Proof:}\\ We prove the lemma indirectly. Then using Hofer's
lemma we can find sequences $\eps_k\searrow 0$ , $z_k\in S$ so
that
\begin{itemize}
\item $\eps_kR_k:=\eps_k|\nabla\tu(z_k)|\longrightarrow+\infty$
\item $|\nabla\tu(z)|\le 2R_k$ whenever $|z-z_k|\le\eps_k$.
\end{itemize}
Writing $z_k=s_k+it_k$, we have to consider the following
situations:
\begin{enumerate}
\item {$-t_kR_k\longrightarrow-\infty$
       \begin{enumerate}
        \item $R_k(1-t_k)\rightarrow l\in[0,+\infty)$
        \item $R-k(1-t_k)\rightarrow+\infty$
       \end{enumerate}}
\item $-t_kR_k\longrightarrow -l\in(-\infty,0]$, then $R_k(1-t_k)\rightarrow+\infty$.
\end{enumerate}
Rescaling in the same way as in the proof of lemma
\ref{lemma1.1.}, i.e.
\[\tu_k(z)=(a(z_k+R_k^{-1}z)-a(z_k),u(z_k+R_k^{-1}z))\ \mbox{for case 1b,}\]
\[\tu_k(z)=(a(s_k+R_k^{-1}z)-a(s_k),u(s_k+R_k^{-1}z))\ \mbox{for case 2}\]
and
\[
\tu_k(z)=(a(s_k+i-R_k^{-1}z)-a(s_k+i),u(s_k+i-R_k^{-1}z))\
\mbox{for case 1a},
\]
we obtain $C^{\infty}_{loc}$--bounds uniform in $k$, where we have
to use the usual elliptic regularity estimates for
$\tu\longmapsto\pas\tu+\tilde{J}(\tu)\pat\tu$ to obtain bounds for
the higher derivatives. Again by the Ascoli--Arzela theorem a
subsequence of $(\tu_k)$ converges to some non constant map
\[
\tilde{w}=(\beta,w):\Omega\longrightarrow\RM,
\]
where $\Omega={\bf C}$ in case 1b and $\Omega=H^+$ in cases 1a and
2. In all these cases we have
\[
\pas\tilde{w}+\tilde{J}(\tilde{w})\pat\tilde{w}=0
\]
and
\[
|\nabla\tilde{w}(z)|\le 2.
\]
In case 2, we have $\tilde{w}(\partial H^+)\subset\RL$, while we
have $\tilde{w}(\partial H^+)\subset \od$ in case 1a. Denote by
$\Omega_k$ the domains of definition of the rescaled maps $\tu_k$,
which are the same as in the proof of lemma \ref{lemma1.1.}. We
claim that
\begin{itemize}
\item $E(\tilde{w})\le E(\tu)$,
\item $\int_{\Omega}w^{\ast}d\lambda=0$.
\end{itemize}
We then have derived a contradiction, because $\tilde{w}$ would
have to be constant (lemma \ref{lemma1.2.} for case 1a,
\cite{A1} for case 2 and \cite{Hofer-Weinstein-conj},\cite{AH} for case 1b).
So let us prove the claim above.\\ Considering case 1b first, we
take $\phi\in\Sigma$ and define $\phi_k\in\Sigma$ by
$$\phi_k(s):=\phi(s-a(z_k)).$$ Then
\begin{eqnarray*}
        \int_{B_{R_k\varepsilon_k}(0)\cap\Omega_k}\tilde{u}_k^{\ast}d(\phi\lambda
        )  & = & \int_{B_{\varepsilon_k}(z_k)\cap({\bf
        R}\times[0,1])}\tilde{u}^{\ast}d(\phi_k\lambda)  \\
         & \le & \int_{{\bf R}\times[0,1]}\tilde{u}^{\ast}d(\phi_k\lambda)  \\
         & \le & E(\tilde{u})
\end{eqnarray*}
Now choose any compact subset $K$ of $\Omega$ and find $k_0\in{\bf
N}$ so that for all $k\ge k_0$ $$K\subset
B_{R_k\varepsilon_k}(0)\cap\Omega_k.$$ Then
$$\int_K\tilde{u}_k^{\ast}d(\phi\lambda)\le E(\tilde{u})\ \forall\
k\ge k_0$$ and therefore $$\int_K\tilde{w}^{\ast}d(\phi\lambda)\le
E(\tilde{u}).$$ Since this holds for all compact subsets $K$ of
$\Omega$ we obtain
$$\int_{\Omega}\tilde{w}^{\ast}d(\phi\lambda)\le E(\tilde{u})$$
and finally taking the supremum over all $\phi\in\Sigma$:
$$E(\tilde{w})\le E(\tilde{u}).$$ Now let $K$ be any compact
subset of $\Omega$. Then for $k$ large enough we have $K\subset
B_{R_k\varepsilon_k}(0)\cap\Omega_k$ and
\begin{eqnarray*}
        \int_K w^{\ast}d\lambda & \le & \left|\int_K w^{\ast}d\lambda-\int_K
        u_k^{\ast}d\lambda\right|+\int_{B_{R_k\varepsilon_k}(0)\cap\Omega_k}
        u_k^{\ast}d\lambda  \\
         & \le & \left|\int_K w^{\ast}d\lambda-\int_K
        v_k^{\ast}d\lambda\right|+\int_{B_{\varepsilon_k}(z_k)\cap({\bf
        R}\times[0,1])}u^{\ast}d\lambda
\end{eqnarray*}
The first term converges to zero for $k\rightarrow+\infty$, but
the second one also does because of $$\int_{{\bf
R}\times[0,1]}u^{\ast}d\lambda=\int_{{\bf
R}\times[0,1]}\tilde{u}^{\ast}d(\phi_0\lambda)\le
E(\tilde{u})<+\infty$$ where $\phi_0\equiv 1\in\Sigma$. This
implies finally $$\int_{\Omega}w^{\ast}d\lambda=0$$ because the
integral vanishes over any compact subset of $\Omega$.\\ In the
cases 2 and 1a the proof of the claim above is essentially the
same up to some minor modifications. We have to define
\[\phi_k(s)=\phi(s-a(s_k))\ \mbox{in case 2}\]
and
\[\Phi_k(s)=\phi(s-a(s_k+i))\ \mbox{in case 1a}.\]
Moreover we have to replace $B_{\eps_kR_k}(0)$ by
$B_{\eps_kR_k}(iR_kt_k)$ in case 2 and
$B_{\eps_kR_k}(iR_k(1-t_k))$ in case 1a respectively. \qed
\begin{proposition}\label{prop5.4.}
Let $\tu$ be a solution of equation (\ref{5.1.1.}). Then every
sequence $(s'_k)_{k\in{\bf N}}\subset{\bf R}$ satisfying
$s'_k\rightarrow +\infty$ or $s'_k\rightarrow-\infty$ has a
subsequence $(s_k)_{k\in{\bf N}}$, so that there is a point
$p\in{\cal L}$ with
\[
\tu(s_k,t)\stackrel{k\rightarrow\infty}{\longrightarrow}(0,p)
\]
in $C^{\infty}([0,1])$.
\end{proposition}
{\bf Proof:}\\ Take any sequence $(s'_k)$ as above and define
\[
\tu_k:S\longrightarrow\RM
\]
by
\[
\tu_k(s,t):=(a(s+s'_k,t)-a(s'_k,0),u(s+s'_k,t)).
\]
Since $\tilde{J}$ does not depend on the ${\bf R}$--component of
$\RM$, we have
\[
\pas\tu_k+\tilde{J}(\tu_k)\pat\tu_k=0.
\]
Moreover with $\tu_k=(a_k,u_k)$:
\[
a_k(0,0)=0,
\]
\[
\tu_k(s,0)\in\RL
\]
and
\[
\tu_k(s,1)\in \od.
\]
Lemma \ref{lemma1.3.} provides a gradient bound for the maps
$\tu_k$ which is uniform in $k$. By elliptic regularity we obtain
uniform $C^{\infty}_{loc}$--bounds and a subsequence of $(\tu_k)$
converges in $C^{\infty}_{loc}$ to some
\[
\tilde{w}=(\beta,w):S\longrightarrow\RM
\]
satisfying
\[
\pas\tilde{w}+\tilde{J}(\tilde{w})\pat\tilde{w}=0,
\]
\[
\tilde{w}(s,0)\in\RL,
\]
\[
\tilde{w}(s,1)\in\od,
\]
\[
\beta(0,0)=0,\]
\[
E(\tilde{w})<+\infty
\]
and
\[
\sup_{z\in S}|\nabla\tilde{w}(z)|<+\infty.
\]
We know that for each $R>0$
\[
\int_{[-R,R]\times[0,1]}u_k^{\ast}d\lambda
\longrightarrow\int_{[-R,R]\times[0,1]}w^{\ast}d\lambda
\]
as $k\rightarrow\infty$. But
\[
\int_{[-R,R]\times[0,1]}u_k^{\ast}d\lambda=\int_{[-R+s_k,R+s_k]\times[0,1]}u^{\ast}d\lambda
\stackrel{k\rightarrow\infty}{\longrightarrow}0,
\]
where $(s_k)$ is a suitable subsequence of $(s'_k)$. This holds
because $u^{\ast}d\lambda$ is a non-negative integrand and
$\int_Su^{\ast}d\lambda\le E(\tu)<+\infty$. Hence
\[
\int_{[-R,R]\times[0,1]}w^{\ast}d\lambda=0
\]
for every $R>0$ and therefore
\[
\int_Sw^{\ast}d\lambda=0.
\]
Lemma \ref{lemma1.1.} implies now that $\tilde{w}$ must be
constant, i.e. $\tilde{w}=(0,w_0)$, where $w_0\in{\cal L}$ might
depend on the sequence $s'_k$ that we chose to define $\tu_k$.
Hence $u(s+s_k,t)\rightarrow w_0$ in $C^{\infty}_{loc}$, in
particular $u(s_k,t)\rightarrow w_0$ in $C^{\infty}([0,1])$.
Moreover $a(s+s_k,t)-a(s_k,0)\rightarrow 0$ in $C^{\infty}_{loc}$.
Choosing $t=1$ we see from the boundary condition $a(s,1)\equiv 0$
that $a(s_k,0)\rightarrow 0$ and therefore $a(s_k,t)\rightarrow 0$
in $C^{\infty}([0,1])$. \qed

\subsection{Existence of an asymptotic limit and exponential decay estimates}
Proposition \ref{prop5.4.} implies that the ends of a finite
energy strip $\tilde{u}$ approach the Legendrian knot $\{0\}\times
{\cal L}\subset {\Real}\times M$. We will go one step further and show that a solution
of equation (\ref{5.1.1.}) has well--defined asymptotic limits. We will also show that the convergence to these asymptotic limits is of exponential nature. The special coordinates derived in proposition \ref{3.1.1.} will
be particularly helpful. 
\begin{proposition}\label{5.1.5.}
Let $\tilde{u}$ be a finite energy strip as in equation
(\ref{5.1.1.}). Then there are points $p_+,p_-\in{\cal L}$ so that
\[
\tilde{u}(s,t)\stackrel{s\rightarrow\pm\infty}{\To} (0,p_{\pm})
\]
in $C^{\infty}([0,1])$.
\end{proposition}
Before we start with the proof of proposition \ref{5.1.5.}, let us
choose convenient coordinates. We will also confine ourselves to
the 'positive end' $s\rightarrow +\infty$ since the negative end
is treated in the same way. By proposition \ref{prop5.4.} we can
find a sequence $s_k\rightarrow +\infty$, so that
$\tilde{u}(s_k,t)$ converges to some point
$(0,p_+)\in\{0\}\times{\cal L}$ in $C^{\infty}([0,1])$ as
$k\rightarrow\infty$ and we may describe $\tilde{u}(s,t)$ by the
coordinates provided by proposition \ref{3.1.1.} if $|s|$ is large
enough. This is because $\tilde{u}(s,t)$ remains near the set $\{0\}\times{\mathcal L}$ for large $|s|$. 
Moreover, our
assumptions imply that the 'ends' of $u$ stay away from the
singular points. We introduce the following change of coordinates
away from the singular points:
\begin{equation}\label{make-boundary-conditions-flat}
{\bf R}\times S^1\times {\bf R}^2\ni
(\tau,\theta,x,y)\longmapsto\left(\tau,\theta,x-\frac{a(\theta)}{b(\theta)}y,y\right).
\end{equation}
We recall (proposition \ref{3.1.1.}) that the spanning surface
${\mathcal D}$ near its boundary is parameterized by
\[
\left\{(\theta,x,y)\in S^1\times{\bf
R}^2\,\left|\,\left(\begin{array}{c}x
\\ y\end{array}\right)=t\cdot\left(\begin{array}{c}a(\theta) \\
b(\theta)\end{array}\right)\,,\,t\in[0,1]\right\}\right.
\]
for suitable functions $a,b:S^1\rightarrow {\bf R}$, and the
singular points correspond to the zeros of $b$.\\ After this
coordinate change we may replace $\Real\times {\cal L}$ by
${\Real}^2\times\{0\}\times\{0\}$, the set $\{0\}\times{\cal
D}^{\ast\ast}$ corresponds to
$\{0\}\times{\Real}\times\{0\}\times{\Real}^{\pm}$ with $\pm=\mbox{sign}(b)$ and we may assume that the point
$(0,p_+)$ corresponds to $0$. Moreover, the contact form $\lambda=dy+xd\theta$ changes to
\begin{equation}\label{new-contact-form}
\hat{\lambda}=dy+(x+\frac{a(\theta)}{b(\theta)}y)d\theta,
\end{equation}
so that the contact structure at the point $(\theta,x,y)$ is generated by
\begin{equation}\label{generators-of-new-contact-str}
\frac{\partial}{\partial\theta}-\left(x+\frac{a(\theta)}{b(\theta)}y\right)\frac{\partial}{\partial y}\ \ \mbox{and}\ \ \frac{\partial}{\partial x}
\end{equation}
and the Reeb vector field changes to
\[
X_{\hat{\lambda}}=\frac{\partial}{\partial y}-\frac{a(\theta)}{b(\theta)}\frac{\partial}{\partial x}.
\]

Our differential equation (\ref{5.1.1.}) has the following form:
\[
v=(\tau,\theta,x,y):[s_0,\infty)\times[0,1]\To{\Real}^4
\]
\begin{equation}\label{5.1.5e.}
\partial_s v+M(v)\partial_t v=0
\end{equation}
\[v(s,0)  \in L_0={\Real}^2\times\{0\}\times\{0\} \]
\[ v(s,1)  \in L_1^=\{0\}\times{\Real}\times\{0\}\times{\Real}.\]
The number $s_0$ is chosen in such a way that $\tilde{u}(s,t)$
lies in the domain in ${\bf R}\times M$ where the above
coordinates exist. The map $M$ is smooth and bounded with values
in $\mbox{GL}({\Real}^4)$, so that all the derivatives are bounded
too and $M^2=-\mbox{Id}$. Because the almost complex structure
$\tilde{J}$ is compatible with the 2--form $\omega$ constructed in
section \ref{3.2.}, we have in addition
\[ M^T\Omega M=\Omega\ \ \mbox{and}\ \ \Omega M>0,\]
where $\Omega$ is a smooth bounded map with bounded derivatives
and values in $\mbox{GL}({\Real}^4)$ so that $\Omega^T=-\Omega$.
We also note that
\begin{equation}\label{5.1.5b.}
\langle v,\Omega(q)\,w\rangle=0
\end{equation}
for $q,v,w\in L_0$ or $q,v,w\in L_1$ since the boundary conditions ${\bf R}\times{\mathcal L}$ and $\{0\}\times{\mathcal D}^{\ast\ast}$ are Lagrangian with respect to the 2--form $\omega$ (here
$\langle\,.\,,\,.\,\rangle$ denotes the standard Euclidean product
on ${\bf R}^4$). Proposition \ref{prop5.4.} implies that
\begin{equation}\label{5.1.5c.}
\sup_{[s,\infty)\times[0,1]}\{|\tau|,|x|,|y|\}\longrightarrow 0
\end{equation}
as $s\rightarrow\infty$, while we only know that
\begin{equation}\label{5.1.5d.}
|\theta(s_k,\,.\,)|_{C^0([0,1])}\To 0
\end{equation}
as $k\rightarrow \infty$. Moreover,
\begin{equation}\label{5.1.6.}
\sup_{[s,\infty)\times[0,1]}|\partial^{\alpha}v|\To 0
\end{equation}
as $s\rightarrow \infty$ for all multi indices $\alpha$ with
$|\alpha|\ge 1$. Our proof of proposition \ref{5.1.5.} consists of
showing that the component $\theta(s,t)$ converges to zero as well
uniformly in $t$, and it will lead also to the following
exponential decay estimates:
\begin{theorem}\label{5.1.5a.}
There exist numbers $\rho,s'>0$ so that we have the following
estimate for each multi index $\alpha\in{\bf N}^2$ , $|\alpha|\ge
0$ and $s\ge s'$:
\[
\sup_{t\in [0,1]}|\partial^{\alpha}v(s,t)|\,\le
c_{\alpha}e^{-\rho(s-s')},
\]
where $c_{\alpha}$ are suitable positive constants.
\end{theorem}
{\bf Proof of proposition \ref{5.1.5.}:}\\
In the following we always assume $s\ge s_0$ so that our boundary
value problem (\ref{5.1.1.}) can be written in coordinates as
(\ref{5.1.5e.}). While we proceed with the proof, it will be
necessary to successively choose a larger constant $s_0$. We will
still denote this constant by $s_0$.\\ We consider the following
family of inner products on $L^2([0,1],{\bf R}^4)$:
\begin{equation}\label{5.1.7.}
(\gamma_1,\gamma_2)_s:=\int^1_0\langle\gamma_1(t),\Omega(v(s,t))
M(v(s,t))\gamma_2(t)\rangle\,dt,
\end{equation}
where $s\ge s_0$ and where $\langle\,.\, ,\, .\,\rangle$ denotes
the Euclidean product on ${\bf R}^4$. We will in future write
$M(s,t)$ and $\Omega(s,t)$ instead of $M(v(s,t))$ and
$\Omega(v(s,t))$. In view of (\ref{5.1.6.}) we have for all
multi indices $\alpha\in{\bf N}^2$ with $|\alpha|\ge 1$
\begin{equation}\label{5.1.8.}
|\partial^{\alpha}\Omega(s,t)|\,
,\,|\partial^{\alpha}M(s,t)|\longrightarrow 0
\end{equation}
uniformly in $t$ as $s$ tends to $+\infty$. Then the norms
$\|\,.\,\|_s$ on $L^2([0,1],{\bf R}^4)$ induced by the products
(\ref{5.1.7.}) are all uniformly equivalent to the usual $L^2$
norm $\|\,.\,\|$, i.e. there are positive constant $c_0,c_1$
independent of $s$ so that
\begin{equation}\label{5.1.7a.}
c_0\|\,.\,\|\le\|\,.\,\|_s\le c_1\|\,.\,\|.
\end{equation}
Define the following dense subspace of $L^2([0,1],{\bf R}^4)$:
\begin{equation}\label{definition-of-H12L}
H_L^{1,2}([0,1],{\bf R}^4):=\{\gamma\in H^{1,2}([0,1],{\bf
R}^4)\ |\ \gamma(0)\in L_0\,,\gamma(1)\in L_1\},
\end{equation} 
where 
\[
L_0:={\bf R}^2\times\{0\}\times\{0\}\ \mbox{and}\ L_1:=\{0\}\times{\bf R}\times\{0\}\times{\bf R}.
\]
In view of the
Sobolev embedding theorem this definition makes sense. We consider
the following family of unbounded linear operators on $L^2$ with
domain of definition $H^{1,2}_L$:
\[ L^2([0,1],{\bf R}^4)\supset H^{1,2}_L([0,1],{\bf
R}^4)\longrightarrow L^2([0,1],{\bf R}^4)\]
\[(A(s)\gamma)(t):=-M(s,t)\dot{\gamma}(t).\]
Since the proof of proposition \ref{5.1.5.} requires some work, we
break up the proof into several lemmas. The following
straightforward lemma summarizes some properties of the operators
$A(s)$:
\begin{lemma}\label{5.1.9.}
The adjoint operator $A(s)^{\ast}$ of $A(s)$ with respect to the
$L^2$--product (\ref{5.1.7.}) has the same domain of definition as
$A(s)$ and is given by
\[(A(s)^{\ast}\gamma)(t)=(A(s)\gamma)(t)-(\Theta(s)\gamma)(t),\]
where $\Theta(s):L^2([0,1],{\bf R}^4)\rightarrow L^2([0,1],{\bf
R}^4)$ is the following zero--order operator:
\[(\Theta(s)\gamma)(t):=M(s,t)\Omega^{-1}(s,t)\partial_t\Omega(s,t)\gamma(t).\]
Moreover, $\Theta(s)(H^{1,2})\subset H^{1,2}$ , $\Theta(s)$ is
antisymmetric and
\begin{equation}\label{5.1.10.}
\|\partial^{k}_s\Theta(s)\|_{{\mathcal L}(L^2,L^2)}\longrightarrow
0\ \mbox{as}\ s\rightarrow\infty,
\end{equation}
where $k\ge 0$.
\end{lemma}\qed \\
Our differential equation (\ref{5.1.1.}) can then be written as
\begin{equation}\label{main-PDE2}
\partial_sv(s,t)\,=\,\big(A(s)v(s)\big)(t),
\end{equation}
with $v(s):=v(s,\,.\,)$. The kernel $\Lambda$ of the operators
$A(s)$ consists of the constant paths with image in $L_0\cap L_1$,
which is a 1--dimensional set. Let
\[P_s:L^2([0,1],{\bf
R}^4)\rightarrow\Lambda\] be the orthogonal projection with
respect to the inner product (\ref{5.1.7.}) and let
\[Q_s:=\,\mbox{Id}\;-P_s.\]
Since the kernels of the operators $A(s)$ all agree, we have the
following important property:
\begin{equation}\label{5.1.11.}
\mbox{The operators}\
\partial_sQ_s\,,\,\partial_{ss}Q_s\ \mbox{have image in}\ \Lambda.
\end{equation}
The following estimate is crucial: 
\begin{lemma}\label{5.1.12.}
There are constants $s_0,\delta>0$ so that for all $s\ge s_0$ and
$\gamma\in H^{1,2}_L([0,1],{\bf R}^4)$ the following inequality
holds:
\[\|A(s)\gamma\|_s\ge\delta\|Q_s\gamma\|_s.\]
\end{lemma}
{\bf Proof:}\\ Proceeding indirectly, we assume that there are
sequences $\delta_k\searrow 0$ , $s_k\nearrow\infty$ and
$\gamma_k\in H^{1,2}_L([0,1],{\bf R}^4)$ so that
\[\|A(s_k)\gamma_k\|_{s_k}<\delta_k\|Q_{s_k}\gamma_k\|_{s_k}.\]
Consider now
\[\alpha_k=\frac{Q_{s_k}\gamma_k}{\|Q_{s_k}\gamma_k\|_{s_k}},\]
so that $0<c_1^{-1}\le\|\alpha_k\|_{L^2}\le c_0^{-1}$ and
\[\|\dot{\alpha}_k\|_{L^2}\le c
\|A(s_k)\alpha_k\|_{s_k}<\delta_k.\] Here we have used that the
norms $\|\,.\,\|_{L^2}$ and $\|\,.\,\|_{s_k}$ are equivalent
(\ref{5.1.7a.}) and that the norm $\|\,.\,\|_{s_k}$ is
$M(s_k)$--invariant. Because the embedding $H^{1,2}([0,1],{\bf
R}^4)\hookrightarrow L^2([0,1],{\bf R}^4)$ is compact, a
subsequence of $(\alpha_k)$ converges in $L^2$ to some $\alpha$.
In view of $\dot{\alpha}_k\stackrel{L^2}{\rightarrow}0$ the
convergence is actually of quality $H^{1,2}$, therefore $\alpha\in
H^{1,2}([0,1],{\bf R}^4)$ and $\dot{\alpha}=0$, i.e.
$\alpha\equiv\,$const. Now $H_L^{1,2}\subset H^{1,2}$ is closed
and $\alpha_k\in H_L^{1,2}([0,1],{\bf R}^4)$, hence
$\alpha\in\Lambda=L_0\cap L_1$.\\ On the other hand, we have $
(\alpha_k,\alpha)_{s_k}=0$ which leads to the contradiction
\begin{eqnarray*}
0<\frac{2}{c_1} & \le & \|\alpha_k\|^2_{L^2}+\|\alpha\|^2_{L^2} \\
 & \le & \frac{1}{c^2_0}(\|\alpha_k\|^2_{s_k}+\|\alpha\|^2_{s_k}) \\
 & = & \frac{1}{c_0^2}(\|\alpha_k-\alpha\|^2_{s_k}) \\
 & \le & \frac{c_1^2}{c_0^2}\|\alpha_k-\alpha\|^2_{L^2} \\
 & \rightarrow & 0
 \end{eqnarray*}
\qed

Let us introduce some notation which will also be useful for
deriving the crucial exponential decay estimate. Fix some integer
$N\ge 1$ and introduce the vector
\[V(s):=(v(s),\pas v(s),\ldots,\partial^{N-1}_s v(s)),\]
which is an element in the N--fold Cartesian product of
$H^{1,2}_L([0,1],{\bf R}^4)$, which we will denote by
$(H^{1,2}_L)^N$. Applying the operator $A(s)$ to each component we
obtain an operator
\[A(s):(H^{1,2}_L)^N\longrightarrow (L^2)^N\]
with
\[\ker A(s)=\Lambda^N.\]
The vector $V$ satisfies the following partial differential
equation:
\begin{equation}\label{5.1.12a.}
\pas V(s)=A(s)V(s)+\hat{\Delta}(s)\pat V(s),
\end{equation}
where
\[\hat{\Delta} (s,t)
  := \left(\begin{array}{cccccc}
      0 & 0 & 0 & \hdots & \hdots & 0 \\
      \Delta_{11}(s,t) & 0 & 0 & \hdots & \hdots & 0 \\
      \Delta_{22}(s,t) & \Delta_{12}(s,t) & 0 & \hdots & \hdots & 0 \\
      \vdots & \vdots & \vdots &  &  & \vdots \\
      \Delta_{N-1,N-1}(s,t) & \Delta_{N-2,N-1}(s,t) &
\Delta_{N-3,N-1}(s,t) & \hdots &
      \Delta_{1,N-1}(s,t) & 0
      \end{array}\right)\]
with
\[\Delta_{lk}(s,t):=\left(\begin{array}{c} k\\l \end{array}\right)
      \partial^l_s ( - M(v(s,t))).\]
The following rather remarkable estimate is essential for the proofs of proposition
\ref{5.1.5.} and theorem \ref{5.1.5a.}. The choices of the inner products in (\ref{5.1.7.}) and lemma \ref{5.1.9.} are crucial for the proof.
\begin{lemma}\label{5.1.13.}
There are numbers $s_0,\delta>0$ so that the function
\[g(s):=\frac{1}{2}\|Q_sV(s)\|^2_s\]
satisfies
\[g''(s)\ge\frac{1}{2}\delta^2\,g(s).\]
\end{lemma}
{\bf Proof:}\\ We have
\[g(s)=\frac{1}{2}\int^1_0\langle
Q_sV(s)(t),\Omega(s,t)M(s,t)Q_sV(s)(t)\rangle\,dt,\] therefore,
using $(\Omega M)^T=\Omega M$,
\[g'(s)=(\pas [Q_sV(s)],Q_sV(s))_s+\frac{1}{2}\int^1_0\langle
Q_sV(s)(t),\pas[\Omega(s,t)M(s,t)]Q_sV(s)(t)\rangle\,dt\] and

\begin{eqnarray*}
g''(s) & = &
\|\pas[Q_sV(s)]\,\|^2_s+(\partial_{ss}[Q_sV(s)],Q_sV(s))_s+
\\
 & & +2\,\int^1_0\langle\pas
 [Q_sV(s)(t)],\pas[\Omega(s,t)M(s,t)]Q_sV(s)(t)\rangle\,dt+\\
  & & +\frac{1}{2}\int^1_0\langle
 Q_sV(s)(t),\partial_{ss}[\Omega(s,t)M(s,t)]Q_sV(s)(t)\rangle\,dt\\
 & =: & T_1+\ldots +T_4 \\
 & \ge & T_2+T_3+T_4.
 \end{eqnarray*}
We can estimate
\begin{equation}\label{5.1.14.}
|T_4|\le \varepsilon(s)\,\|Q_sV(s)\|^2_s,
\end{equation}
where $0<\varepsilon(s)\stackrel{s\rightarrow\infty}{\longrightarrow}
0$ is a suitable function. From now on, we will write
$\varepsilon(s)$ for any positive function which decays to zero as
$s\rightarrow\infty$.\\ Now let us estimate $T_3$. We have to
consider $\pas P_s$ first. If $e\in L_0\cap L_1$ then $P_s\gamma$
is given by
\begin{equation}\label{5.1.28.}
P_s\gamma=\frac{(\gamma,e)_s}{\|e\|^2_s}\cdot e
\end{equation}
and
\begin{eqnarray}\label{5.1.15.}
|(\pas P_s)\gamma| & = & \left|\frac{\frac{\partial}{\partial s
}(\gamma,e)_s}{\|e\|^2_s} -
2\frac{(\gamma,e)_s}{\|e\|^3_s}\frac{\partial}{\partial s
}\|e\|_s\cdot e\right| \nonumber \\ & \le &
\varepsilon(s)\|\gamma\|_s,
\end{eqnarray}
since there are positive constants $c_0,c_1$ so that
$c_0\le\|e\|_s\le c_1$ for all $s$. Moreover,
\begin{eqnarray}\label{5.1.16.}
|\pas[\Omega(s)M(s)]\gamma| & = & |[D\Omega(v(s))\pas
v(s)]M(s)\gamma+\Omega(v(s))DM(v(s))[\pas v(s),\gamma]| \nonumber
\\ & \le & c\,|\pas v(s)|\,|\gamma|,
\end{eqnarray}
where $c>0$ is some constant. Using (\ref{5.1.15.}),
(\ref{5.1.16.}),
\begin{itemize}
\item
\[\|Q_s\pas
V(s)\|_{C^0([0,1])}\stackrel{s\rightarrow\infty}{\longrightarrow}0,\]
which follows from (\ref{5.1.6.}),
\item $\|V(s)\|_{L^2}$ is bounded uniformly in $s$ and
\item $\pas Q_s+\pas P_s=0$,
\end{itemize}
we obtain
\begin{eqnarray}\label{5.1.17.}
|T_3| & = & \left|\int^1_0\langle Q_s\pas V(s)(t)-[\pas P_s]
V(s)(t),\pas[\Omega(s,t)M(s,t)]Q_sV(s)(t)\rangle\,dt\right|
\nonumber
\\
 & \le & c\|Q_s\pas V(s)-[\pas P_s] V(s)\|_{C^0([0,1])}\|\pas
 v(s)\|_{L^2}\|Q_sV(s)\|_{L^2}  \nonumber \\
 & = & \varepsilon(s)\|A(s)v(s)\|_{L^2}\|Q_sV(s)\|_{L^2} \\
 & \le & \varepsilon(s)\|A(s)V(s)\|_{L^2}\|Q_sV(s)\|_{L^2}.\nonumber
\end{eqnarray}
We are now left with $T_2$. Shortening the notation, we write
$\pas Q_s\gamma$ instead of $(\pas Q_s)\gamma$ and
$Q_s\hat{\Delta}(s)\pat V(s)$ instead of $Q_s(\hat{\Delta}(s)\pat
V(s))$ etc. We calculate
\begin{eqnarray*}
\pas(Q_sV(s)) & = & \pas
Q_s\,V(s)+Q_sA(s)V(s)+Q_s\hat{\Delta}(s)\pat V(s) \\
 & = & \pas
 Q_s\,V(s)+A(s)Q_sV(s)-P_sA(s)V(s)+Q_s\hat{\Delta}(s)\pat V(s)
\end{eqnarray*}
and
\begin{eqnarray*}
\partial_{ss}(Q_sV(s)) & = & \partial_{ss}Q_s\,V(s)+\pas Q_s\,\pas
V(s)+A(s)\pas V(s) - \pas M(s)\,\pat V(s)+\\
 & & +\pas Q_s\hat{\Delta}(s)\pat V(s)
 +Q_s\pas\hat{\Delta}(s)\,\pat
 V(s)+Q_s\hat{\Delta}(s)\partial_{st}V(s)-\\
  & & -\pas P_sA(s)V(s)+P_s\pas M(s)\,\pat V(s)-P_sA(s)\pas V(s).
\end{eqnarray*}
We write the term $\hat{\Delta}(s)\partial_{st}V(s)$ as
$\tilde{\Delta}(s)\pat V(s)$, where
\[
\tilde{\Delta} (s,t):=\left(
    \begin{array}{ccccc}
      0        & 0              & 0               & \dots  & 0 \\
      0        & \Delta_{11} (s,t)    & 0               & \dots & 0 \\
      0        & \Delta_{22} (s,t)    & \Delta_{12}(s,t)      & \dots & 0 \\
      \vdots   & \vdots         & \vdots          & \ddots & \vdots \\
      0        & \Delta_{N-1,N-1}(s,t) & \Delta_{N-2,N-1}(s,t)& \dots &
\Delta_{1,N-1}(s,t)
    \end{array}\right).
\]
Inserting this into $T_2$ we obtain with (\ref{5.1.11.})
\begin{eqnarray*}
T_2 & = & (A(s)\pas V(s),Q_sV(s))_s-(\pas M(s)\,\pat
V(s),Q_sV(s))_s+ \\
 & & +(Q_s\pas\hat{\Delta}(s)\pat
V(s),Q_sV(s))_s+(Q_s\tilde{\Delta}(s)\pat V(s),Q_sV(s))_s \\
 & =: & T_{21}+T_{22}+T_{23}+T_{24}.
 \end{eqnarray*}
 We estimate
 \[|T_{22}|=|(\pas
 M(s)\,M(s)\,A(s)Q_sV(s),Q_sV(s))_s|\le\varepsilon(s)\|A(s)Q_sV(s)\|_s\|Q_sV(s)\|_s.\]
 The expressions $T_{23},T_{24}$ are estimated similarly, so that
 \begin{equation}\label{5.1.18.}
 |T_{22}|,|T_{23}|,|T_{24}|\le\varepsilon(s)\|A(s)Q_sV(s)\|_s\|Q_sV(s)\|_s.
\end{equation}
Using $\pat V(s)=M(s)A(s)Q_sV(s)$, equation (\ref{5.1.12a.}) and lemma \ref{5.1.9.}, we continue with $T_{21}$:
\begin{eqnarray*}
T_{21} & = & (\pas V(s),A(s)Q_sV(s))_s-(\pas
V(s),\Theta(s)Q_sV(s))_s \\
 & = & \|A(s)Q_sV(s)\|^2_s+(\hat{\Delta}(s)\pat
 V(s),A(s)Q_sV(s))_s -\\
  & & -(A(s)Q_sV(s),\Theta(s)Q_sV(s))_s-(\hat{\Delta}(s)\pat
  V(s),\Theta(s)Q_sV(s))_s\\
  & \ge &
  \|A(s)Q_sV(s)\|^2_s-\varepsilon(s)\|A(s)Q_sV(s)\|^2_s-\varepsilon(s)\|A(s)Q_sV(s)\|_s
 \|Q_sV(s)\|_s \\
 & \ge & \frac{1}{2} \|A(s)Q_sV(s)\|^2_s-\varepsilon(s)\|A(s)Q_sV(s)\|_s
 \|Q_sV(s)\|_s\ \mbox{for large }s.
\end{eqnarray*}
Using lemma \ref{5.1.12.}, the inequalities (\ref{5.1.14.}),
(\ref{5.1.17.}), (\ref{5.1.18.}) and the above estimate for
$T_{21}$, we obtain
\begin{eqnarray*}
g''(s) & \ge & T_{21}-|T_{22}|-|T_{23}|-|T_{24}|-|T_3|-|T_4|
\\
 & \ge & \frac{1}{2} \|A(s)Q_sV(s)\|^2_s-\varepsilon(s)\|A(s)Q_sV(s)\|_s
 \|Q_sV(s)\|_s-\varepsilon(s)\|Q_sV(s)\|^2_s \\
& = & \|A(s)Q_sV(s)\|_s\left(\frac{1}{2}
\|A(s)Q_sV(s)\|_s-\varepsilon(s)\|Q_sV(s)\|_s\right)-\varepsilon(s)\|Q_sV(s)\|^2_s
\\
& \ge &
\|A(s)Q_sV(s)\|_s\|Q_sV(s)\|_s\left(\frac{\delta}{2}-\varepsilon(s)\right)-\varepsilon(s)\|Q_sV(s)\|^2_s
\\
& \ge &
\left(\frac{\delta^2}{3}-\varepsilon(s)\right)\|Q_sV(s)\|_s^2\ \
\mbox{, where }s\mbox{ is so large that
$\frac{\delta}{2}-\varepsilon(s)\ge\frac{\delta}{3}$} \\ & \ge &
\frac{\delta^2}{4}\|Q_sV(s)\|^2_s\ \ \ s\mbox{ so large that
$\frac{\delta^2}{3}-\varepsilon(s)\ge\frac{\delta^2}{4}$} \\ & = &
\frac{\delta^2}{2}g(s).
\end{eqnarray*}
This completes the proof of lemma \ref{5.1.13.}.\qed

\begin{lemma}\label{5.1.26.}
Let $s_0,\delta$ be as in lemma \ref{5.1.13.}. Then we have for
all $s\ge s_1\ge s_0$
\[
g(s)\le g(s_1)e^{-\frac{\delta}{\sqrt{2}}(s-s_1)}.
\]
\end{lemma}
{\bf Proof:}\\ Defining
$h(s):=g(s)-g(s_1)e^{-\frac{\delta}{\sqrt{2}}(s-s_1)}$, we observe
that $h(s_1)=0$ and $h''(s)\ge\frac{\delta^2}{2}h(s)$ in view of
lemma \ref{5.1.13.}, hence $h$ cannot have a local maximum with
$h>0$. On the other hand, we also have $h(s)\rightarrow 0$ as
$s\rightarrow\infty$ in view of $g''(s)\rightarrow 0$ and lemma
\ref{5.1.13.}. We conclude $h\le 0$ which proves the lemma. \qed
\\

We now have to estimate $|P_sv(s)|$ and $|P_s\pas v(s)|$, the
components of $v(s)$ and $\pas v(s)$ along $\Lambda=\ker A(s)$.
\begin{lemma}\label{5.1.27.}
If $s\ge s_0$ then
\[
|P_s\pas v(s)|\le \varepsilon(s)\,\|Q_sv(s)\|_s,
\]
where $0<\varepsilon(s)\rightarrow 0$ as $s\rightarrow \infty$.
\end{lemma}
{\bf Proof:}\\ We compute using (\ref{5.1.7a.}),
(\ref{main-PDE2}), (\ref{5.1.28.}) and lemma \ref{5.1.9.}
\begin{eqnarray*}
|P_s\pas v(s)| & \le & c\,|(\pas v(s),e)_s| \\
 & = & c\,|(A(s)v(s),e)_s| \\
 & = & c\,|(A(s)Q_sv(s),e)_s| \\
 & = & c\,|(Q_sv(s),-\Theta(s)e)_s| \\
 & \le & \varepsilon(s)\,\|Q_sv(s)\|_s
\end{eqnarray*}
\qed \\

{\bf Proof of proposition \ref{5.1.5.} (continued):}\\ We will
show now that
\[
|P_sv(s)|\longrightarrow 0
\]
as $s\rightarrow\infty$. Then we are done because
\[
\|v(s)\|_{L^2}\le
|P_sv(s)|+\|Q_sv(s)\|_{L^2}\stackrel{s\rightarrow
\infty}{\longrightarrow} 0
\]
and
\[
\|\pas v(s)\|_{L^2}\stackrel{s\rightarrow \infty}{\longrightarrow}
0,
\]
i.e. $v(s)$ converges to zero in $H^{1,2}([0,1])$ and therefore
also in $C^0([0,1])$ by the Sobolev embedding theorem.\\

In view of equations (\ref{5.1.7a.}),(\ref{5.1.28.}) we have to show that
\[
|(v(s),e)_s|\longrightarrow 0
\]
as $s\rightarrow\infty$. We know already that
\[
|(v(s_k),e)_{s_k}|\le\|v(s_k)\|_{s_k}\|e\|_{s_k}\le c
\,\|v(s_k)\|_{L^2}\stackrel{k\rightarrow\infty}{\longrightarrow}
0.
\]
We estimate for $s\ge s_k$, combining lemma \ref{5.1.26.} and
lemma \ref{5.1.27.}, with $c$ being a generic constant independent
of $k$ and $s$
\begin{eqnarray*}
|(v(s),e)_s-(v(s_k),e)_{s_k}| & = &
\left|\int_{s_k}^s\frac{d}{d\sigma}(v(\sigma),e)_{\sigma}d\sigma\right|
\\
& \le &
c\int_{s_k}^s\|\partial_{\sigma}v(\sigma)\|_{\sigma}d\sigma+\\ & &
+\int_{s_k}^s\int_0^1|\langle
v(\sigma,t),\partial_{\sigma}[\,\Omega(v(\sigma,t))M(v(\sigma,t))\,]\cdot
e\rangle |\,dt\,d\sigma\\ & \le &
c\,\|\partial_sv(s_k)\|_{s_k}\int_{s_k}^se^{-\frac{\delta}{2\sqrt{2}}(\sigma-s_k)}\,d\sigma+\\
 & & +c\,\int_{s_k}^s\|v(\sigma)\|_{L^2([0,1])}\|\partial_{\sigma}v(\sigma)\|_{L^2([0,1])}d\sigma\\
& \le & c\|\pas
v(s_k)\|_{s_k}(1-e^{-\frac{\delta}{2\sqrt{2}}(s-s_k)})\\ & \le &
c\|\pas v(s_k)\|_{s_k},
\end{eqnarray*}
which converges to zero if $k\rightarrow\infty$. This completes
the proof of proposition \ref{5.1.5.}. \qed
\\

{\bf Proof of theorem \ref{5.1.5a.}:}\\ We saw earlier that lemmas
\ref{5.1.26.} and \ref{5.1.27.} imply
\[
\|\pas v(s)\|_{L^2}\le c\, e^{-\frac{\delta}{2\sqrt{2}}(s-s_0)}
\]
for all $s\ge s_0$, where $c,s_0>0$ are suitable constants. In
view of $\pas v(s,t)+M(v(s,t))\pat v(s,t)=0$ we also have
\[
\|\pat v(s)\|_{L^2}\le c\,e^{-\frac{\delta}{2\sqrt{2}}(s-s_0)}
\]
for a suitable positive constant $c$. Note that
$v(s)=-\int^{+\infty}_s\pas v(\sigma)\,d\sigma$ so that
\begin{eqnarray*}
\|v(s)\|_{L^2} & \le & \int^{+\infty}_s\|\pas
v(\sigma)\|_{L^2}d\sigma \\ & \le & c\,\int^{+\infty}_s
e^{-\frac{\delta}{2\sqrt{2}}(\sigma-s_0)}d\sigma \\ & = &
\frac{2c\sqrt{2}}{\delta}e^{-\frac{\delta}{2\sqrt{2}}(s-s_0)}.
\end{eqnarray*}
Hence we know already that $\|\partial^{\alpha}v(s)\|_{L^2}$
decays exponentially fast with rate at least
$\rho=-\frac{\delta}{2\sqrt{2}}$ whenever $|\alpha|\le 1$. Because
of the Sobolev embedding theorem we obtain exponential decay for
$\sup_{0\le t\le 1}|v(s)|$ as well. We have to use induction to
obtain the same decay behavior for the higher derivatives of $v$.
Recalling that we defined
\[
V(s)=(v(s),\pas v(s),\ldots,\partial_s^{N-1}v(s)), \ N\ge 1,
\]
we know that $\|Q_s\pas^k v(s)\|_{L^2}$ exhibits the desired
exponential decay for any integer $k$. Assume that
$\|V(s)\|_{L^2}$ decays exponentially with rate $\rho$ as above
(we know that this is true for $N=2$). We claim that then $\|\pas
V(s)\|_{L^2}$ and $\|\pat V(s)\|_{L^2}$ have to decay
exponentially with the same rate as well. Applying $Q_s$ to
equation (\ref{5.1.12a.}) and multiplying with $M(s)$ we obtain
\[
\pat V(s)=M(s)Q_s\pas
V(s)+M(s)P_sA(s)V(s)-M(s)Q_s\hat{\Delta}(s)\pat V(s),
\]
which implies
\[
\|\pat V(s)\|_{L^2}\le c\,\|Q_s\pas
V(s)\|_{L^2}+c\,|P_sA(s)V(s)|+\varepsilon(s)\|\pat V(s)\|_{L^2},
\]
i.e. for $s$ so large that $\varepsilon(s)\le 1/2$
\begin{equation}\label{5.1.29.}
\|\pat V(s)\|_{L^2}\le 2c\,\|Q_s\pas
V(s)\|_{L^2}+2c\,|P_sA(s)V(s)|.
\end{equation}
The expression $\|Q_s\pas V(s)\|_{L^2}$ decays exponentially by
lemma \ref{5.1.26.} and the other also does because of
\begin{eqnarray*}
|P_sA(s)V(s)| & \le & c\,|(A(s)V(s),e)_s| \\
 & \le & c\,|(V(s),\Theta(s)e)_s| \\
 & \le & \varepsilon(s)\,\|V(s)\|_{L^2},
\end{eqnarray*}
where $e\in \Lambda^N$. This proves our claim, i.e. we have now
shown exponential decay for $\| \pat^k\pas^l v(s)\|_{L^2}$, where
$k\in\{0,1\}$ and $l\ge 0$ is an arbitrary integer.\\ Equation
(\ref{5.1.12a.}) yields
\[
\pat V(s)=(\mbox{Id}+M(s)\hat{\Delta}(s))^{-1}M(s)\,\pas V(s)
\]
(the inverse makes sense if $s$ is sufficiently large), and
differentiating the above identity successively by $t$ shows by
induction that $\|\pat^kV(s)\|_{L^2}$ decays exponentially for
arbitrary integers $k$. The desired decay for the $C^0$ norm then
follows from the Sobolev embedding theorem.\qed

\subsection{An asymptotic formula}

We need to know more about the asymptotic behavior of the solutions than merely the apriori estimate in theorem \ref{5.1.5a.}. The aim is to
prove the asymptotic formula (theorem \ref{asymptotic-formula-theorem}):

\begin{theorem}
For sufficiently large $s_0$ and $s\ge s_0$ we have the following
asymptotic formula for non constant solutions $v$ of
(\ref{5.1.5e.}) having finite energy:
\begin{equation}\label{asymptotic-formula}
v(s,t)=e^{\int_{s_0}^s\alpha(\tau)d\tau}\Big(e(t)+r(s,t)\Big),
\end{equation}
where $\alpha:[s_0,\infty)\rightarrow{\bf R}$ is a smooth function
satisfying $\alpha(s)\rightarrow\lambda < 0$ as
$s\rightarrow\infty$ with $\lambda$ being an eigenvalue of the
selfadjoint operator
\[
A_{\infty}:L^2([0,1],{\bf R}^4)\supset H^{1,2}_L([0,1],{\bf
R}^4)\longrightarrow L^2([0,1],{\bf R}^4)
\]
\[
\gamma\longmapsto -M_{\infty}\dot{\gamma}\ ,\
M_{\infty}:=\lim_{s\rightarrow\infty}M(v(s,t))
\]
(see (\ref{definition-of-H12L}) for the definition of the domain of $A_{\infty}$).
Moreover, $e(t)$ is an eigenvector of $A_{\infty}$ belonging to
the eigenvalue $\lambda$ with $e(t)\neq 0$ for all $t\in[0,1]$, 
and $r$ is a smooth function so that $r$
and all its derivatives converge to zero uniformly in $t$ as
$s\rightarrow\infty$.
\end{theorem}
{\bf Remark:}\\
The above theorem is of course also valid for the negative end, $s\rightarrow-\infty$, of a solution. We have the same formula as in (\ref{asymptotic-formula}), but the function $\alpha(s)$ will converge to a positive eigenvalue of the operator $A_{\infty}=-M_{\infty}\frac{d}{dt}$, where $M_{\infty}=\lim_{s\rightarrow-\infty}M(v(s,t))$.\\

The first step in the proof is the following proposition. The steps from proposition \ref{5.1.19.} below to theorem
\ref{asymptotic-formula-theorem} are very similarly to the corresponding results in \cite{A1} or \cite{HWZI}.

\begin{proposition}\label{5.1.19.}
There is a number $s_0>0$ so that
\[
\|v(s)\|_s=e^{\int_{s_0}^s\alpha(\tau)d\tau}\|v(s_0)\|_{s_0}
\]
for all $s\ge s_0$, where $\alpha$ has the properties stated above
in theorem \ref{asymptotic-formula-theorem}.
\end{proposition}

Before we can continue with the proof, we need some information
about the spectra of the selfadjoint operators
$A(s)-\frac{1}{2}\Theta(s)$.

%
%
%
\begin{theorem}\label{spectral-gaps}
For each $L>0$ there are numbers $d,s_1>0$ and a sequence
$r_n\in[nL,(n+1)L]\ ,\ n\in{\bf Z}$ so that
\[
[r_n-d,r_n+d]\cap\sigma(A(s)-\frac{1}{2}\Theta(s))=\emptyset
\]
for all $s\ge s_1$.
\end{theorem}
{\bf Proof:}\\ Let us review the strategy of the proof: We want to view $A(s)-\frac{1}{2}\Theta(s)$ as a perturbation of $A_{\infty}$, the operator obtained for $s\rightarrow\infty$. There are theorems about the spectrum of selfadjoint operators in a Hilbert space perturbed by bounded symmetric operators. The trouble here is that $A_{\infty}-A(s)+\frac{1}{2}\Theta(s)$ is not a bounded operator. We fix this by introducing operators $B_{\infty}$ and $B(s)$, all having the same first order term, and which are unitary equivalent to the operators $A_{\infty}$ and $A(s)$ so that it suffices to study the spectra of $B_{\infty}$ and $B(s)$. \\

We would like to find a smooth map
\[
T:[s_0,\infty)\times[0,1]\longrightarrow\mbox{GL}({\bf R}^4)
\]
so that $T(s,\,.\,)$ converges in $C^{\infty}([0,1])$ to some
$T_{\infty}\in\mbox{GL}({\bf R}^4)$ satisfying the following
conditions:
\begin{itemize}
\item $T^tT=\Omega M$,
\item $TM=J_0T$,
\item $T^tJ_0T=-\Omega$,
\end{itemize}
with corresponding conditions for $T_{\infty}$ as
$s\rightarrow\infty$. Here $T^t$ denotes the transpose of $T$ and
$J_0$ is multiplication by $i$ on ${\bf C}^2$ if we identify ${\bf
R}^4$ with ${\bf C}^2$. Actually two of the above conditions imply
the third one. We may view the map $T$ as a unitary trivialization
of the hermitian vector bundle
\[
\Big(([s_0,\infty)\times[0,1])\times{\bf
R}^4,\Omega,M\Big)\widetilde{\rightarrow}\Big(([s_0,\infty)\times[0,1])\times{\bf
R}^4,-J_0,J_0\Big).
\]
The construction of $T$ is Gram--Schmidt orthogonalization with respect to the hermitian bundle metric
\[
h=\langle\,.\,,\,\Omega M\,.\,\rangle+i\langle\,.\,,\,\Omega\,.\,\rangle.
\]
We define $T(s,t)$ by mapping the generator
\[
\sigma(s,t):=\frac{\partial}{\partial\theta}-\left(x(s,t)+\frac{a(\theta(s,t))}{b(\theta(s,t))}y(s,t)\right)\frac{\partial}{\partial y}
\]
of the contact structure (\ref{generators-of-new-contact-str}) in $u(s,t)$ onto $(0,1)\in{\bf C}^2$. Consequently, the maps $$ \Phi_s:(L^2([0,1],{\bf
R}^4),(.,.)_s)\rightarrow(L^2([0,1],{\bf R}^4) ,(.,.)_{L^2}) $$ $$
\gamma\longmapsto T(s,.)\gamma $$ $$ \Phi_{\infty}:(L^2([0,1],{\bf
R}^4),(.,.)_s|_{s\rightarrow\infty})\rightarrow(L^2([0,1],{\bf
R}^4),(.,.)_{L^2}) $$ $$ \gamma\longmapsto T_{\infty}\gamma $$ are
isometries. They map $H^{1,2}_{L}([0,1],{\bf R}^4)$ onto $$
H^{1,2}_{L_s}([0,1],{\bf R}^4):=\left\{\gamma\in
H^{1,2}([0,1],{\bf R}^4) \left|\begin{array}{c} \gamma(0)\in
T(s,0)\cdot L_0
\\ \gamma(1) \in T(s,1)\cdot L_1
\end{array}\right.\right\}
$$ and $$ H^{1,2}_{L_{\infty}}([0,1],{\bf R}^4):=\left\{\gamma\in
H^{1,2}([0,1],{\bf R}^4) \left|\begin{array}{c} \gamma(0)\in
T_{\infty}\cdot L_0 \\ \gamma(1) \in T_{\infty}\cdot L_1
\end{array}\right.\right\} $$ respectively. We consider the
following operators $$ \tilde{B}(s):L^2([0,1],{\bf R}^4)\supset
H^{1,2}_{L_{s}}([0,1],{\bf R}^4) \rightarrow L^2([0,1],{\bf R}^4)
$$ $$ \tilde{B}(s):=\Phi_s\circ
(A(s)-\frac{1}{2}\Theta(s))\circ\Phi_s^{-1} $$ $$
B_{\infty}:L^2([0,1],{\bf R}^4)\supset
H^{1,2}_{L_{\infty}}([0,1],{\bf R}^4) \rightarrow L^2([0,1],{\bf
R}^4) $$ $$ B_{\infty}:=\Phi_{\infty}\circ
A_{\infty}\circ\Phi_{\infty}^{-1} $$ where we equip
$L^2([0,1],{\bf R}^4)$ with the ordinary $L^2$-inner product
$(.,.)_{L^2}$. Unitary equivalent selfadjoint operators have the
same spectrum, hence $$
\sigma(\tilde{B}(s))=\sigma(A(s)-\frac{1}{2}\Theta(s)) $$ and $$
\sigma(B_{\infty})=\sigma(A_{\infty}). $$ It remains to
investigate the spectra of $\tilde{B}(s)$ and $B_{\infty}$. First
we note that the operators $\tilde{B}(s)$ and $B_{\infty}$ are
selfadjoint with respect to the standard $L^2$--product. Let us
compute them. We obtain
\begin{equation}
\tilde{B}(s)  =  -J_0\frac{\partial}{\partial t}+J_0\frac{\partial
T(s)}{\partial t}T(s)^{-1}-\frac{1}{2}T(s)\Theta(s) T(s)^{-1},
\end{equation}
where the operator $-J_0\frac{\partial}{\partial t}$ is
selfadjoint and the operator
\[
S(s):L^2([0,1],{\bf R}^4)\longrightarrow L^2([0,1],{\bf R}^4)
\]
\[
\gamma\longmapsto J_0\frac{\partial T(s)}{\partial
t}T(s)^{-1}\gamma-\frac{1}{2}T(s)\Theta(s) T(s)^{-1}\gamma
\]
is symmetric. We note that $S(s)$ converges to zero as
$s\rightarrow\infty$ in the operator norm. The operator
$B_{\infty}$ is simply given by $-J_0\frac{\partial}{\partial
t}$.\\

Summarizing, we have introduced coordinates so that the operators
$A(s)-\frac{1}{2}\Theta(s)$ and $A_{\infty}$ correspond to
operators with the same first order term on the same Hilbert space $(L^2([0,1],{\bf R}^n), (\,.\, ,\,.\,)_{L^2})$, but they all have
different domains of definition. We have to fix this without
changing anything that we have achieved so far.\\

We can find a smooth map
\[
C:[s_0,\infty)\times[0,1]\longrightarrow
\mbox{Sp}(4)\cap\mbox{O}(4)=\mbox{U}(2)
\]
having the following properties:
\begin{itemize}
\item $C(s,\,.\,)\rightarrow\,\mbox{Id}$ in $C^{\infty}([0,1])$ as
$s\rightarrow\infty$,
\item $C(s,0)T(s,0)L_0=T_{\infty}L_0$,
\item $C(s,1)T(s,1)L_1=T_{\infty}L_1$.
\end{itemize}
The operators
\[
B(s):L^2([0,1],{\bf R}^4)\supset H^{1,2}_{L_{\infty}}([0,1],{\bf
R}^4)\longrightarrow L^2([0,1],{\bf R}^4)
\]
\[
(B(s)\gamma)(t):=C(s,t)^{-1}(\tilde{B}(s)C(s)\gamma)(t)
\]
have the form
\[
B(s)=B_{\infty}+\Delta(s,t),
\]
where $\gamma\longmapsto\Delta(s)\gamma$ is a symmetric zero order
perturbation with $\|\Delta(s)\|\rightarrow 0$ as
$s\rightarrow\infty$ in the operator norm. They are unitary
equivalent to $\tilde{B}(s)$ hence the spectra are the same. The spectrum of the operator $B_{\infty}$, which has
domain of definition $H^{1,2}_{L_{\infty}}([0,1],{\bf R}^4)$ consists of all integer multiples of $\pi/2$. Moreover, the spectrum consists of eigenvalues only since the resolvent of $B_{\infty}$ is a compact operator. Every eigenvalue has multiplicity one. Verifying this is a straight forward computation which we leave to the reader. Let us summarize our discussion as follows:
\begin{proposition}\label{spectrum-of-A-infty}
The spectrum of the operator
\[
A_{\infty}:L^2([0,1],{\bf R}^4)\supset H^{1,2}_L([0,1],{\bf
R}^4)\longrightarrow L^2([0,1],{\bf R}^4)
\]
\[
\gamma\longmapsto -M_{\infty}\dot{\gamma}\ ,\
M_{\infty}:=\lim_{s\rightarrow\infty}M(v(s,t))
\]
consists of all integer multiples of $\frac{\pi}{2}$. The resolvent of the operator $A_{\infty}$ is a compact operator on $L^2([0,1],{\bf R}^4)$. All the points in the spectrum are eigenvalues of multiplicity one.
\end{proposition}
\qed\\

In order to control the spectra of the perturbations $B(s)$ we will need the following perturbation result (see \cite{A1}) which follows
from a result of T. Kato (see \cite{Kato}):
\begin{theorem}\label{Kato}
Let $T:H\supset D(T)\rightarrow H$ be a selfadjoint operator in a
Hilbert space $H$ and let $A_0:H\rightarrow H$ be a linear, bounded
and symmetric operator. Then the following holds:
\begin{itemize}
    \item  \begin{eqnarray*}
             & & \mbox{dist}(\sigma(T),\sigma(T+A_0)) \\
             & := & \max\big\{\sup_{\lambda\in\sigma(T)}\mbox{
             dist}(\lambda, \sigma(T+A_0))\, , \,
             \sup_{\lambda\in\sigma(T+A_0)}\mbox{
             dist}(\lambda,\sigma(T))\big\} \\
             & \le & \|A_0\|_{{\cal L}(H)}
            \end{eqnarray*}

    \item  Assume further that the resolvent $(T-\lambda_0)^{-1}$ of $T$
    exists and is compact for some $\lambda_0\not\in \sigma(T)$.\\
    Then $(T-\lambda)^{-1}$ is compact for every $\lambda\not\in\sigma(T)$
    and $\sigma(T)$ consists of isolated eigenvalues $\{\mu_k\}_{k\in{\bf
    Z}}$ with finite multiplicities $\{m_k\}_{k\in{\bf Z}}$.\\
    If we assume that $\sup_{k\in{\bf Z}}m_k\le M<\infty$ and that for each
    $L>0$ there is a number $m_T(L)\in{\bf N}$ so that every interval $I\subset
    {\bf R}$ of length $L$ contains at most $m_T(L)$ points of $\sigma(T)$
    (counted with multiplicity) then for each $L>0$ there is also a number
    $m_{T+A_0}(L)\in{\bf N}$ so that every interval $I\subset
    {\bf R}$ of length $L$ contains at most $m_{T+A_0}(L)$ points of $\sigma(T+A_0)$.
\end{itemize}
\end{theorem}\qed

We find for all $L > 0$ some $m \in {\bf N}$ so that every
interval $I \subseteq {\bf R}$ of length $L$ contains at most $m$
points of the spectrum of $B_{\infty}$.

Moreover by theorem \ref{Kato},
\begin{equation} \label{*}
  \mbox{ dist } (\sigma (B_{\infty}), \sigma (B(s)) \rightarrow 0
\end{equation}
as $s \rightarrow \infty$.

Define now the intervals $$
  I_n := [ n L,  (n+1) L]\ ; \quad n \in {\bf Z}.
$$ Then each $I_n$ contains at most $m$ points of $\sigma
(B_{\infty})$, so there is a closed subinterval $J_n \subset I_n$
of length $\frac{L}{m+1}$ that does not contain any point of
$\sigma (B_{\infty})$. Because of (\ref{*}) there is a closed
interval $J'_n \subseteq J_n \subseteq I_n$ of length $\frac{L}{2
(m+1)}$ which does not contain any point of $\sigma (B(s))$
whenever $s \geq s_1$ where $s_1$ is sufficiently large (this
$s_1$ does not depend on $n$).

So we found a sequence $r_n \in I_n$ and a positive constant $d$,
so that $$
  [r_n - d, r_n + d] \cap \sigma (B(s)) = \emptyset
$$ for all large $s$. This completes the proof of theorem
\ref{spectral-gaps}.\qed
%
%
%
%

{\bf Proof of proposition \ref{5.1.19.}:}\\ This result has an analogue
in \cite{A1} and \cite{HWZI},\cite{HWZIV}. However, there are some
different features due to the boundary condition and the
degeneracy of the problem. We assume first that
$\|v(s,\,.\,)\|_{C^0([0,1])}\neq 0$ if $s$ is sufficiently large.
As in the references cited above, it is very easy to state the
correct function $\alpha$ so that we have the proposed formula for
$\|v(s)\|_s$. Indeed, we have to take
\[
\alpha(s):=\frac{\frac{d}{ds}\|v(s)\|^2_s}{2\|v(s)\|^2_s}.
\]
We define now
\[
\xi(s,t):=\frac{v(s,t)}{\|v(s)\|_s}
\]
and note that
\begin{equation}\label{5.1.20.}
\pas\xi(s,t)+M(s,t)\pat\xi(s,t)+\alpha(s)\xi(s,t)=0.
\end{equation}
We define
\[\Gamma_1(s,t):=-\frac{1}{2}M(s,t)\Omega^{-1}(s,t)\pas(\Omega
M)(s,t)\] and the covariant derivative
\[
\nabla_s\xi(s):=\partial_s\xi(s)+\Gamma_1(s)\xi(s)
\]
so that for all smooth $u_1,u_2:{\bf R}\rightarrow L^2([0,1],{\bf
R}^4)$
\[
\frac{d}{ds}(u_1(s),u_2(s))_s=(\nabla_s
u_1(s),u_2(s))_s+(u_1(s),\nabla_s u_2(s))_s,
\]
hence
\begin{equation}\label{5.1.20a.}
0=(\nabla_s\xi(s),\xi(s))_s.
\end{equation}
The partial differential equation for $\xi$ can be written in the
form
\begin{equation}\label{5.1.20b.}
A(s)\xi(s)=\nabla_s\xi(s)+\alpha(s)\xi(s)-\Gamma_1(s)\xi(s)
\end{equation}
which implies
\begin{equation}\label{5.1.21.}
\alpha(s) = (\xi(s),\Gamma_1(s)\xi(s))_s+(\xi(s),A(s)\xi(s))_s.
\end{equation}
We define
\[
 \Gamma_2(s,t):=-M(s,t)\Omega^{-1}(s,t)\pas\Omega(s,t)
\,M(s,t).
\]
and
\[
\Gamma_3(s,t):=\Omega^{-1}(s,t)\pas\Omega(s,t).
\]
Computing the adjoint operators $\Gamma_1^{\ast}$ and
$\Gamma_2^{\ast}$ with respect to the inner product (\ref{5.1.7.})
yields
\[\Gamma_1^{\ast}=\Gamma_1\ \mbox{and}\
\Gamma_2^{\ast}=\Gamma_3.\] Introducing the operator
\begin{eqnarray*}
\Gamma_4(s)\xi(s) & := & (\nabla_sM(s))M(s)\xi(s) \\
 & := & -M(s)\nabla_s(M(s)\xi(s))-\nabla_s\xi(s) \\
 & = & \frac{1}{2}(\Gamma_3-\Gamma_2),
\end{eqnarray*}
we find that
\[
\Gamma_4^{\ast}(s)=\frac{1}{2}(\Gamma_3^{\ast}-\Gamma_2^{\ast})=-\Gamma_4.
\]
A simple calculation shows also that
\begin{equation}\label{5.1.20c.}
\pat\nabla_s-\nabla_s\pat=\pat\Gamma_1.
\end{equation}
Using now the partial differential equation (\ref{5.1.20b.}) ,
lemma \ref{5.1.9.}, equations
(\ref{5.1.10.}),(\ref{5.1.20a.}),(\ref{5.1.20c.}) and the fact
that $\|\Gamma_k(s)\xi(s)\|^2_s\rightarrow 0\ ,\ k=1,\ldots,4\,,$
as $s\rightarrow\infty$ we estimate the derivative of $\alpha$ as
follows:
\begin{eqnarray*}
\alpha'(s) & = &
(\nabla_s(A(s)\xi(s)),\xi(s))_s+(A(s)\xi(s),\nabla_s\xi(s))_s+\\
 & &
 (\nabla_s(\Gamma_1(s)\xi(s)),\xi(s))_s+(\Gamma_1(s)\xi(s),\nabla_s\xi(s))_s
 \\
 & =: & T_1+T_2+T_3+T_4.
\end{eqnarray*}
We have
\[
|T_4|\le\varepsilon(s)\|\nabla_s\xi(s)\|_s
\]
and
\begin{eqnarray*}
|T_3| & \le &
|((\nabla_s\Gamma_1(s))\xi(s),\xi(s))_s|+|(\Gamma_1(s)\nabla_s\xi(s),\xi(s))_s|
\\
& \le & \varepsilon(s)+\varepsilon(s)\|\nabla_s\xi(s)\|_s.
\end{eqnarray*}
Inserting (\ref{5.1.20b.}) and using (\ref{5.1.20a.}), we obtain
\begin{eqnarray*}
T_2 & = &
\|\nabla_s\xi(s)\|^2_s-(\Gamma_1(s)\xi(s),\nabla_s\xi(s))_s\\
 & \ge &
 \|\nabla_s\xi(s)\|^2_s-\varepsilon(s)\|\nabla_s\xi(s)\|_s.
\end{eqnarray*}
We now take care of the term $T_1$:
\begin{eqnarray*}
T_1 & = &
-(M(s)\nabla_s(M(s)A(s)\xi(s)),\xi(s))_s-(\Gamma_4(s)A(s)\xi(s),\xi(s))_s
\\
& = &
(-M(s)\nabla_s\pat\xi(s),\xi(s))_s+(A(s)\xi(s),\Gamma_4(s)\xi(s))_s
\\
& = &
(A(s)\nabla_s\xi(s),\xi(s))_s+(M(s)\pat\Gamma_1(s)\xi(s),\xi(s))_s+\\
& & +(A(s)\xi(s),\Gamma_4(s)\xi(s))_s \\ & = &
(\nabla_s\xi(s),A(s)\xi(s))_s-(\nabla_s\xi(s),\Theta(s)\xi(s))_s+
\\
& &
+(M(s)\pat\Gamma_1(s)\xi(s),\xi(s))_s+(\nabla_s\xi(s),\Gamma_4(s)\xi(s))_s+
\\
& &
+\alpha(s)(\xi(s),\Gamma_4(s)\xi(s))_s-(\Gamma_1(s)\xi(s),\Gamma_4(s)\xi(s))_s
\\
& =:& T_{11}+\ldots+T_{16}.
\end{eqnarray*}
The term $T_{11}$ is identical with $T_2$ which we estimated
above. The expressions $|T_{12}|$, $|T_{14}|$ can be estimated
from above by $\varepsilon(s)\|\nabla_s\xi(s)\|_s$ while
$|T_{13}|$ and $|T_{16}|$ tend to zero as $s\rightarrow\infty$.
The term $T_{15}$ vanishes since $\Gamma_4$ is skew--adjoint.
Summarizing, we got the following inequality for the derivative of
$\alpha$:
\begin{eqnarray}\label{5.1.22.}
\alpha'(s) & \ge & T_1+T_2-|T_3|-|T_4|  \nonumber \\
  & \ge & 2\|\nabla_s\xi(s)\|^2_s-\varepsilon(s)\|\nabla_s\xi(s)\|_s-\varepsilon(s).
\end{eqnarray}

We assume now that the function $\alpha$ is not bounded from above
and we wish to derive a contradiction. Then we can find a sequence
$s_k\rightarrow \infty$ so that $\alpha(s_k)\rightarrow\infty$. If
we had $\alpha(s)\ge\eta>0$ for all large $s$ and some positive
number $\eta$ then we would obtain
\[\|v(s)\|_{L^2}\ge c\|v(s)\|_s\ge
e^{\eta(s-s_0)}\|v(s_0)\|_{s_0}\rightarrow\infty\] in
contradiction to the fact that $|v(s,\,.\,)|\rightarrow 0$
uniformly in $t$ as $s\rightarrow\infty$. Because of theorem
\ref{spectral-gaps} we may pick $\eta>0$ in such a way that there
is a positive number $d$ so that $\eta-d>0$ and
\[
[\eta-d,\eta+d]\cap\sigma(A(s)-\frac{1}{2}\Theta(s))=\emptyset.
\]
Then we can find a sequence $s'_k\rightarrow\infty$ so that
$\alpha(s'_k)<\eta$. We may also assume that
$s'_k<s_{k+1}<s'_{k+1}$ and $\alpha(s_k)>\eta$. Hence, if $\alpha$
is not bounded from above then it must oscillate. Let $\hat{s}_k$
be the smallest number with $\hat{s}_k>s_k$ and
$\alpha(\hat{s}_k)=\eta$. Since the operators
$A(s)-\frac{1}{2}\Theta(s)$ are selfadjoint we have for every
$\theta$ in the resolvent set
\begin{equation}\label{norm-of-the-resolvent}
\|(A(s)-\frac{1}{2}\Theta(s)-\theta\,\mbox{Id})^{-1}\|_s=
\frac{1}{\mbox{dist}(\theta,\sigma(A(s)-\frac{1}{2}\Theta(s)))}.
\end{equation}
Recalling the differential equation (\ref{5.1.20b.}) for $\xi$, we
obtain ($\varepsilon_k$ being a suitable sequence of positive
numbers converging to zero)
\begin{eqnarray*}
1 & = & \|\xi(\hat{s}_k)\|_{\hat{s}_k} \\
 & \le &
 \|(A(\hat{s}_k)-\frac{1}{2}\Theta(\hat{s}_k)-\eta\,\mbox{Id})^{-1}\|_{\hat{s}_k}\|\nabla_s\xi(\hat{s}_k)
 -\Gamma_1(\hat{s}_k)\xi(\hat{s}_k)-\frac{1}{2}\Theta(\hat{s}_k)\xi(\hat{s}_k)\|_{\hat{s}_k}
  \\
 & \le & \frac{1}{d}\|\nabla_s\xi(\hat{s}_k)
 -\Gamma_1(\hat{s}_k)\xi(\hat{s}_k)-\frac{1}{2}\Theta(\hat{s}_k)\xi(\hat{s}_k)\|_{\hat{s}_k}
  \\
 & \le &
 \frac{1}{d}\|\nabla_s\xi(\hat{s}_k)\|_{\hat{s}_k}+\varepsilon_k,
\end{eqnarray*}
i.e. for sufficiently large $k$
\begin{equation}\label{s-deriv-bd-from-below}
0<\frac{d}{2}\le\|\nabla_s\xi(\hat{s}_k)\|_{\hat{s}_k}.
\end{equation}
We now insert this into inequality (\ref{5.1.22.}) and obtain that
for sufficiently large $k$
\[
\alpha'(\hat{s}_k)>0,
\]
which would imply $\alpha(s)<\eta$ for $s<\hat{s}_k$ close to
$\hat{s}_k$ in contradiction to the definition of $\hat{s}_k$.
Hence $\alpha$ must be bounded from above.\\

Let us show now that $\alpha$ cannot be unbounded from below either.
Pick a sequence $r_n$ as in theorem \ref{spectral-gaps}. Assuming
in the contrary that $\alpha$ is not bounded from below we can
find $s_n$ so that $\alpha(s_n)=r_n$ and $\alpha'(s_n)<0$. In the
same way as we derived (\ref{s-deriv-bd-from-below}), we also
obtain here
\[
0<\frac{d}{2}\le\|\nabla_s\xi(s_n)\|_{s_n}
\]
for all large $n$ and
\[
\alpha'(s_n)>0
\]
which is a contradiction. Therefore $\alpha$ must also be bounded
from below.\\

There exists a sequence $s_k\rightarrow\infty$ so that
$\|\nabla_s\xi(s_k)\|_{s_k}\rightarrow 0$. Otherwise we had
$\|\nabla_s\xi(s)\|_s\ge\eta>0$ for a suitable $\eta$ and all
large $s$ which would imply $\alpha'(s)\ge\frac{1}{2}\eta^2$ for
all large $s$ and $\alpha(s)\rightarrow\infty$ as
$s\rightarrow\infty$ which is not true.\\
Because $\alpha$ is bounded, we can find a subsequence (which we also
denote by $(s_k)_{k \in {\bf N}}$) so that
$$
  \lim_{k \rightarrow \infty} \alpha (s_k) = \lambda
$$
exists. We claim that $\lambda \in \sigma (A_{\infty})$.
If we had $\lambda \not\in \sigma (A_{\infty})$ then $\varepsilon :=
\inf_{\mu \in \sigma (A_{\infty})} |{\lambda - \mu} |> 0$ because
$\sigma (A_{\infty})$ is closed and therefore
$$
  |{\mu' - \lambda} |\geq \varepsilon - |{\mu - \mu'}| \ \ \forall \mu \in
  \sigma (A_{\infty}), \mu' \in \sigma (A(s) -\frac{1}{2}\Theta(s))
$$
which implies
$$
  \mbox{ dist } (\lambda, \sigma (A(s)-\frac{1}{2}\Theta(s))) \geq \varepsilon -
  \sup_{\mu' \in \sigma
  (A(s)-\frac{1}{2}\Theta(s))} \mbox{ dist } (\mu' , \sigma (A_{\infty})) > \varepsilon/2
$$
if $s$ is sufficiently large, by theorem \ref{spectral-gaps}, i.e.
$$
 \alpha (s_k) \not\in \sigma (A(s_k) -\frac{1}{2}\Theta(s_k))
$$
for $k$ sufficiently large.

Then
\begin{eqnarray*}
1 & = & \|\xi(s_k)\|_{s_k} \\
 & = & \|(A(s_k)-\frac{1}{2}\Theta(s_k)-\alpha(s_k)\,\mbox{Id})^{-1}(\nabla_s\xi(s_k)-\frac{1}{2}\Theta(s_k) -\Gamma_1(s_k)\xi(s_k))\|_{s_k} \\
 & \le & \frac{4}{\varepsilon}\|\nabla_s\xi(s_k)\|_{s_k}+\varepsilon_k,
\end{eqnarray*}
where $k$ is chosen so large that $|{\lambda - \alpha (s_k)}| <
\varepsilon/4$ and $\varepsilon_k\searrow 0$ is a suitable sequence.
But this contradicts $\|\nabla_s \xi (s_k) \|_{s_k} \rightarrow
0$, hence $\lambda \in \sigma (A_{\infty})$.\\

Let us show that indeed
$$
  \lim_{s \rightarrow \infty} \alpha (s) = \lambda .
$$
Take now a sequence $s_k \rightarrow \infty$ and assume that
there are subsequences $(s'_k),(s''_k)$ which converge to
different limits $\lambda'$ and $\lambda''$. By our previous
discussion we have
\[
\lambda',\lambda''\in\sigma(A_{\infty})
\]
and we assume that $\lambda'<\lambda''$. We may also assume that
$s'_k<s''_k<s'_{k+1}$. It is a consequence of theorem
\ref{spectral-gaps} that there are $d>0$ and
$\nu\in(\lambda',\lambda'')$ so that
\[
\mbox{dist}(\nu,\sigma(A(s)-\frac{1}{2}\Theta(s)))\ge d
\]
whenever $s$ is sufficiently large. Let now $s$ be any number with
$\alpha(s)=\nu$. Then we estimate as before:
\begin{eqnarray*}
1 & = & \|\xi(s)\|_{s} \\
 & = & \|(A(s)-\frac{1}{2}\Theta(s)-\nu\,\mbox{Id})^{-1}(\nabla_s\xi(s)-\frac{1}{2}\Theta(s) -\Gamma_1(s)\xi(s))\|_{s} \\
 & \le & \frac{1}{d}\|\nabla_s\xi(s)\|_{s}+\varepsilon(s),
\end{eqnarray*}
where $\varepsilon(s)$ is a suitable positive function tending to
zero as $s\rightarrow\infty$. Using inequality (\ref{5.1.22.}), we
obtain $\alpha'(s)>0$ for all large enough $s$ with
$\alpha(s)=\nu$, but this is a contradiction since it prohibits
$\alpha$ from oscillating between $\lambda'$ and $\lambda''$.
Hence the limit
\[
\lambda=\lim_{s\rightarrow\infty}\alpha(s)\in\sigma(A_{\infty})
\]
exists and it is indeed an eigenvalue because the operator $A_{\infty}$
has compact resolvent so that the spectrum consists of eigenvalues
only. Moreover, $\lambda\le 0$ since otherwise
$\|v(s)\|_{L^2}\rightarrow\infty$. Let us show that $\lambda<0$\\

We know that there are $\rho,s_0>0$ so that for all $s\ge s_0$:
\[
\|v(s)\|_s\le ce^{-\rho(s-s_0)}
\]
which follows from theorem \ref{5.1.5a.}. Using proposition \ref{5.1.19.}, we see that the function 
\[
e^{\rho(s-s_0)}\|v(s)\|_s = \|v(s_0)\|_{s_0}e^{\rho(s-s_0)+\int_{s_0}^s\alpha(\tau)d\tau}
\]
remains bounded for all $s\ge s_0$. This means that the function
\[
f(s):=\rho(s-s_0)+\int_{s_0}^s\alpha(\tau)d\tau
\]
has to be bounded as well. Now
\[
f'(s)=\rho+\alpha(s)\rightarrow \rho+\lambda
\]
as $s\rightarrow\infty$. Boundedness of $f$ implies then $\rho+\lambda\le 0$.\\

It remains to take care of the case for which $\|v(s)\|_s=0$ for
some $s$. Then $v(s,t)=0$ for all $t\in[0,1]$ and a simple
application of the similarity principle implies that $v$ is
constant (see \cite{A1}, \cite{HZbook}) in contradiction to our
assumptions. This completes the proof of proposition
\ref{5.1.19.}.\qed

The following three lemmas are versions of lemmas in \cite{A1} and
\cite{HWZI}. The proof of theorem \ref{asymptotic-formula-theorem}
is then very similar to the corresponding version in \cite{A1}.
For the convenience of the reader, we sketch the path until the
proof of theorem \ref{asymptotic-formula-theorem}. The proofs of
the corresponding results in \cite{A1} and \cite{HWZI} can almost
be carried over verbatim; we will indicate the necessary
modifications.

\begin{lemma}\label{lemma A}
For every $\beta = (\beta_1, \beta_2) \in {\bf N}^2$ and $j \in
{\bf N}$ we have
\[
  \sup_{(s,t) \in [s_0, \infty) \times [0,1]}
    |\partial^{\beta} \xi (s,t)| < \infty
    \]
    \[
  \sup_{s_0 \leq s < \infty}
    \left|\frac{d^j \alpha}{ds^j} (s)\right| < \infty
\]
where $\xi (s,t) = \frac{v(s,t)}{\|v(s)\|_s}$ and $\alpha (s) =
(A(s) \cdot \xi (s) + \Gamma_1 (s) \cdot \xi (s), \xi (s))_s$.\\
(here, we adopt the convention $0\in{\bf N}$ )
\end{lemma}
{\bf Proof:}\\ This is actually a version of lemma 3.10. from
\cite{A1}. The proof remains essentially the same. There are two
minor modifications: The operator $T_{\infty}(t)$ in \cite{A1}
should be replaced by the $t$--independent operator $T_{\infty}$
that we introduced in the proof of theorem \ref{spectral-gaps}.
Moreover, the estimate for $|\alpha'(s)|$ in \cite{A1} has to be
replaced by
\[
|\alpha'(s)|\le
c'\|\pas\xi(s)\|^2_{L^2([0,1])}+c''\|\pas\xi(s)\|_{L^2([0,1])}+c''',
\]
which follows from the estimates that lead us to inequality
(\ref{5.1.22.}). We then get for $p>2$ and $\delta_2>0$
\begin{eqnarray*}
 \|\alpha'\|^p_{L^p([s^{\ast}-\delta_2,s^{\ast}+\delta_2])}
 & \le &
 4^{p-1}(c')^p\int_{s^{\ast}-\delta_2}^{s^{\ast}+\delta_2}\left(\int_0^1|\pas\xi(s,t)|^2dt\right)^p
 ds+ \\
  & & +4^{p-1}(c'')^p\int_{s^{\ast}-\delta_2}^{s^{\ast}+\delta_2}\left(\int_0^1|\pas\xi(s,t)|^2\right)^{p/2}
 ds+\\
 & & +2\cdot 4^{p-1}(c''')^p\delta_2 \\
 & \le &
 4^{p-1}(c')^p\|\pas\xi\|^{2p}_{L^{2p}(Q_{\delta_2})}+4^{p-1}(c'')^p\|\pas\xi\|^{p}_{L^{2p}(Q_{\delta_2})}+\\
 & & +2\cdot 4^{p-1}(c''')^p\delta_2,
\end{eqnarray*}
where
$Q_{\delta_2}:=[s^{\ast}-\delta_2,s^{\ast}+\delta_2]\times[0,1]$;
but this estimate works as well as the original one in
\cite{A1}.\qed

\begin{lemma}\label{lemma B}
Let $$
  E \subseteq H^{1,2}_{L} ([0,1], {\bf R}^4)
   \subseteq L^2 ([0,1], {\bf R}^4)
$$ be the eigenspace of $A_{\infty}$ belonging to $\lambda \in
\sigma (A_{\infty})$.

Then $$
  \inf_{e \in E} \|\xi (s) - e\|_{H^{1,2} ([0,1], {\bf R}^4)}
  \rightarrow 0
$$ as $s \rightarrow \infty$.
\end{lemma}

{\bf Proof:}\\ This is a modification of lemma 3.6. in
\cite{HWZI}. The proof is very similar to \cite{HWZI}, replace
$\pas\xi$ in the estimates by the covariant derivative
$\nabla_s\xi$.\qed

\begin{lemma}\label{lemma C}
There exists $e\in E$ such that $\xi(s)\rightarrow e$ in
$H^{1,2}([0,1],{\bf R}^4)$ as $s\rightarrow\infty$.
\end{lemma}
{\bf Proof:}\\ This is essentially lemma 3.12. in \cite{A1}. Using
the $L^2$--product
\[
(u_1,u_2):=\int^1_0\langle
u_1(t),\Omega_{\infty}M_{\infty}u_2(t)\rangle\,dt
\]
instead, the proof in \cite{A1} can be carried over.\qed

{\bf Proof of theorem \ref{asymptotic-formula-theorem}:}\\ By
proposition \ref{5.1.19.} we have
\begin{eqnarray*}
v(s,t) & = & \|v(s)\|_s\xi(s,t) \\
 & = & e^{\int_{s_0}^s\alpha(\tau)d\tau}\|v(s_0)\|_{s_0}\xi(s,t)
 \\
 & = & e^{\int_{s_0}^s\alpha(\tau)d\tau}[\tilde{e}(t)+r(s,t)]
\end{eqnarray*}
with
\[
r(s,t):=\|v(s_0)\|_{s_0}(\xi(s,t)-e(t)),
\]
\[
\tilde{e}(t):=\|v(s_0)\|_{s_0}e(t)\in E,
\]
where $e(t)$ is the eigenvector given by lemma \ref{lemma C}. Recall from the proof of theorem \ref{spectral-gaps} that the operator $A_{\infty}$ is unitary equivalent to the operator $B_{\infty}=-i\frac{d}{dt}$ acting on a suitable closed subspace of $H^{1,2}([0,1],{\bf R}^4)$. Eigenvectors of $B_{\infty}$ are of the form $\hat{e}(t)=e^{i\lambda t}\hat{e}(0)$, hence they are nowhere zero and so are eigenvectors of $A_{\infty}$. The
proof that $r$ converges to zero in $C^{\infty}$ is the same as in
\cite{A1}, so we omit the details.\qed

\subsection{Proof of theorem \ref{convergence-of-alpha}}

We will need later the following simple observation concerning the function $\alpha$ which appears in the asymptotic formula, theorem \ref{asymptotic-formula-theorem}:
\begin{proposition}\label{alphaprime-goes-to-zero}
All derivatives of the function $\alpha$ as in (\ref{5.1.5e.}) converge to zero as $|s|\rightarrow \infty$.
\end{proposition}

{\bf Proof:}\\
We have $\|\pas^k\xi(s)\|_{L^2([0,1])}\rightarrow 0$ for $k\ge 1$ and $s\rightarrow\infty$ because $\pas^k\xi(s,t)$ equals up to multiplication with a constant the derivative $\pas^k r(s,t)$, where $r$ is the remainder in the asymptotic formula, theorem \ref{asymptotic-formula-theorem}. 
Recall equation (\ref{5.1.21.}):
\begin{eqnarray*}
\alpha(s) & = & (\xi(s),\Gamma_1(s)\xi(s))_s+(\xi(s),A(s)\xi(s))_s
\end{eqnarray*}
Differentiating with respect to $s$, we obtain the assertion of the proposition.
\qed \\

We denote by $E$ the eigenspace of the asymptotic operator $A_{\infty}$ belonging to the eigenvalue $\lambda$. Let $e$ be the generator of $E$ such that $\xi(s)\rightarrow e$ as $s\rightarrow\infty$ (see lemma \ref{lemma C}). Let
\[
\pi_s:(L^2([0,1],{\bf R}^4),(\,.\,,\,.\,)_s)\longrightarrow E
\]
\[
\pi_sv:=\frac{(v,e)_s}{\|e\|^2_s}\cdot e
\]
be the orthogonal projection onto the space $E$ and let $Q_s:=\,\mbox{Id}-\pi_s$. The following lemma is similar to lemma \ref{5.1.12.}
\begin{lemma}\label{5.1.12A.}
There are constants $s_0,\delta>0$ so that for all $s\ge s_0$ and
$\gamma\in H^{1,2}_L([0,1],{\bf R}^4)$ the following inequality
holds:
\[\|(A(s)-\alpha(s))Q_s\gamma\|_s\ge\delta\|Q_s\gamma\|_s.\]
\end{lemma}
{\bf Proof:}\\ Proceeding indirectly, we assume that there are
sequences $\delta_k\searrow 0$ , $s_k\nearrow\infty$ and
$\gamma_k\in H^{1,2}_L([0,1],{\bf R}^4)$ so that
\[
\frac{\|(A(s_k)-\alpha(s_k))Q_{s_k}\gamma_k\|_{s_k}}{\|Q_{s_k}\gamma_k\|_{s_k}}<\delta_k.
\]
With
\[
\eta_k:=\frac{Q_{s_k}\gamma_k}{\|Q_{s_k}\gamma_k\|_{s_k}}
\]
we get $\delta_1\ge\|\eta_k\|_{L^2([0,1])}\ge\delta_0>0$ for some $\delta_0,\delta_1$ and
\begin{equation}\label{5.1.12A-1}
\|(A(s_k)-\alpha(s_k))\eta_k\|_{s_k}<\delta_k\rightarrow 0
\end{equation}
We estimate 
\[
\|\pat\eta_k\|_{L^2([0,1])}\le c\,\|A(s_k)\eta_k\|_{s_k}<c(\delta_k+|\alpha(s_k)|)\le 2c|\lambda|,
\]
for sufficiently large $k$. The sequence $\eta_k$ is therefore bounded in $H^{1,2}([0,1],{\bf R}^4)$ which is compactly embedded into $L^2([0,1],{\bf R}^4)$. Hence we may assume that after passing to a suitable subsequence $\eta_k\rightarrow\eta$ in $L^2$. We estimate
\begin{eqnarray*}
\|(A_{\infty}-\lambda)\eta_k\|_{s_k} & \le & \|(A(s_k)-\alpha(s_k))\eta_k\|_{s_k}+\|(A_{\infty}-A(s_k))\eta_k\|_{s_k}+\\
 & & +c\,\|(\alpha(s_k)-\lambda)\eta_k\|_{s_k}\\
 & \le & \delta_k+\|M_{\infty}-M(s_k)\|_{L^{\infty}([0,1])}\|\pat\eta_k\|_{s_k}+c\,|\alpha(s_k)-\lambda|\\
 & \rightarrow & 0
\end{eqnarray*}
and
\begin{eqnarray*}
\|-M_{\infty}\pat\eta_k-\lambda\eta\|_{L^2([0,1])}\le c\,\|(A_{\infty}-\lambda)\eta_k\|_{s_k}+|\lambda|\|\eta-\eta_k\|_{L^2([0,1])},
\end{eqnarray*}
which converges to zero, hence $\pat\eta_k$ converges in $L^2$ to $\lambda M_{\infty}\eta$ which is then the weak derivative $\pat\eta$ of $\eta$. We conclude $A_{\infty}\eta=\lambda\eta$, i.e. $Q_{s_k}\eta=0$ for all $k$. This leads to the contradiction 
\[
1=(\eta_k,\eta_k)_{s_k}\le |(\eta_k,\eta_k-\eta)_{s_k}|+|(\eta_k,\eta)_{s_k}|\le\,\mbox{const}\,\|\eta_k-\eta\|_{L^2([0,1])}\rightarrow 0
\]
and completes the proof of the lemma.
\qed \\

Our aim is now to estimate $\xi(s,t)-e(t)$, $\alpha(s)-\lambda$ and all its derivatives in absolute value from above by $c\,e^{-\delta\,s}$. For an integer $N\ge 1$ we introduce the vector
\[
V(s,t):=\Big(\pas^{k}(\xi(s,t)-e(t))\Big)_{0\le k\le N-1}
\]
and we want derive a PDE satisfied by $V$. Using equation (\ref{5.1.20.}), which is 
\[
\pas\xi(s)=A(s)\xi(s)-\alpha(s)\xi(s),
\]
and $A_{\infty}e=\lambda e$ we obtain
\[
\pas(\xi(s)-e)=(A(s)-A_{\infty})(\xi(s)-e)+(A(s)-A_{\infty})e+(\lambda-\alpha(s))e.
\]
Differentiating successively with respect to $s$ and viewing $\pi_s$, $Q_s$, $A(s)-\alpha(s)$ as operators on $N$--tuples in $(H^{1,2}_L)^N$ we obtain the following differential equation for $V$:
\begin{equation}\label{PDE-for-V}
\pas V(s)=(A(s)-\alpha(s))V(s)-\tilde{\alpha}(s)\,V(s)+H(s)+E(s),
\end{equation} 
where $H$ and its derivatives decay like $e^{-|\lambda|s}$, the vector $E(s)$ is given by
\[
E(s)=\Big(\pas^k(\alpha(s)-\lambda)\cdot e\Big)_{0\le k\le N-1}
\]
and
\[
\tilde{\alpha}(s)=\left(\begin{array}{ccccc} 0 & 0 & 0 & \cdots & 0 \\ \alpha_{11} & 0 & 0 & \cdots & 0 \\ \alpha_{22} & \alpha_{12} & 0 & \cdots & 0 \\ \vdots & \vdots & \vdots & & \vdots \\ \alpha_{N-1,N-1} & \alpha_{N-2,N-1} & \alpha_{N-3,N-1} & \cdots & 0 \end{array}\right)
\]
with
\[
\alpha_{lk}=\left(\begin{array}{c} k \\ l \end{array}\right)\frac{d^l\alpha}{ds^l}\ ,\ 1\le l,k\le N-1.
\]
We also note that
\[
\pi_s E(s)=E(s).
\]
We define now the function
\[
g(s):=\frac{1}{2}\|\pi_sV(s)-V(s)\|^2_s
\]
and we denote by $\Gamma(s)$ a matrix whose entries are zero order operators such that $\|D^{\alpha}\Gamma(s)\|\le c\,e^{-|\lambda|s}$ in the operator norm. We will always use this notation if we are not concerned with the explicit structure of $\Gamma(s)$. We compute
\begin{eqnarray*}
g'(s) & = & (\pas(Q_sV(s)),Q_sV(s))_s+(Q_sV(s),\Gamma(s)Q_sV(s))_s.
\end{eqnarray*}
We continue with the second derivative
\begin{eqnarray*}
g''(s) & = & (\partial_{ss}(Q_sV(s)),Q_sV(s))_s+\|Q_sV(s)\|^2_s+\\
 & & +(Q_sV(s),\Gamma(s)Q_sV(s))_s+(\pas(Q_sV(s)),\Gamma(s)Q_sV(s))_s\\
& \ge & (Q_s(\partial_{ss}V(s)),Q_sV(s))_s+(\pas(Q_sV(s)),\Gamma(s)Q_sV(s))_s+\\
 & & +(Q_sV(s),\Gamma(s)Q_sV(s))_s\\
 & =: & T_1+T_2+T_3,
\end{eqnarray*}
and we note that
\begin{equation}\label{the-t23-estimates} 
|T_2|,|T_3|\le c\,e^{-|\lambda|s}\|Q_sV(s)\|_s.
\end{equation}
We note that $\Gamma(s)$ here is different than in the equation for $g'(s)$, but we use the same symbol since we only care about the exponential decay. We have also used that $V$ and its derivatives are bounded by lemma \ref{lemma A} and that the operators $\pas Q_s$, $\partial_{ss}Q_s$ have range in $E^N$, hence the ranges of these operators are orthogonal to the range of $Q_s$. We will also use the facts that $\pas Q_s-Q_s\pas$ has range in $E^N$ and that $\tilde{\alpha}(s)Q_s-Q_s\tilde{\alpha}(s)=0$. Differentiating (\ref{PDE-for-V}) yields
\begin{eqnarray*}
\partial_{ss}V(s) & = & -\pas(M(v))M(v)A(s)V(s)-\alpha'(s)V(s)+(A(s)-\alpha(s))\pas V(s)-\\
 & & -\tilde{\alpha}'(s)V(s)-\tilde{\alpha}(s)\pas V(s)+H'(s)+E'(s).
\end{eqnarray*}
We evaluate
\begin{eqnarray*}
T_1 & = & (-\pas(M(v))M(v)A(s)V(s),Q_sV(s))_s+((A(s)-\alpha(s))\pas V(s),Q_sV(s))_s-\\
 & & -(\alpha'(s)Q_sV(s),Q_sV(s))_s-(\tilde{\alpha}'(s)Q_sV(s),Q_sV(s))_s+\\
 & & +(H'(s),Q_sV(s))_s\\
& = & T_{11}+\cdots T_{15}.
\end{eqnarray*}
If $\ve(s)$ denotes a function which converges to zero with all derivatives as $s\rightarrow\infty$ then we can estimate
\[
|T_{13}|,|T_{14}|\le\ve(s)\|Q_sV(s)\|^2_s
\]
and
\[
|T_{11}|,|T_{15}|\le c\,e^{-|\lambda|s}\|Q_sV(s)\|_s.
\]
We continue with the term $T_{12}$:
\begin{eqnarray}\label{evaluate-T12}
T_{12}  & =  & ((A(s)-\alpha(s))\pi_s(\pas V(s)),Q_sV(s))_s+\nonumber \\
 & & +((A(s)-\alpha(s))Q_s(\pas V(s)),Q_sV(s))_s.
\end{eqnarray}
We compute
\begin{eqnarray*}
(A(s)-\alpha(s))\pi_s(\pas V(s)) & = & [-M(v)\pat-\alpha(s)]\frac{(\pas V(s),e)_s}{\|e\|^2_s}\cdot e\\
 & = & \frac{(\pas V(s),e)_s}{\|e\|^2_s}(A(s)-\alpha(s))\,e\\
 & = & \frac{(\pas V(s),e)_s}{\|e\|^2_s}\Big((A(s)-A_{\infty})e+(\lambda-\alpha(s))e\Big),
\end{eqnarray*}
so that
\begin{eqnarray}\label{evaluate-T12-1}
& & \Big|((A(s)-\alpha(s))\pi_s(\pas V(s)),Q_sV(s))_s\Big|\nonumber\\
 & = & \left|\frac{(\pas V(s),e)_s}{\|e\|^2_s}((A(s)-A_{\infty})e,Q_sV(s))_s\right|\nonumber\\
 & \le & c\,e^{-|\lambda|s}\|Q_sV(s)\|_s.
\end{eqnarray}
In a similar fashion, we obtain
\begin{equation}\label{evaluate-T12-2}
\Big|((A(s)-\alpha(s))\pi_s V(s),Q_sV(s))_s\Big| \le  c\,e^{-|\lambda|s}\|Q_sV(s)\|_s.
\end{equation}
We evaluate now the second term in equation (\ref{evaluate-T12}) using the differential equation (\ref{PDE-for-V}), equation (\ref{evaluate-T12-2}) and proposition \ref{alphaprime-goes-to-zero}:
\begin{eqnarray*}
& & ((A(s)-\alpha(s))Q_s(\pas V(s)),Q_sV(s))_s\\ & = & (Q_s(\pas V(s)),(A(s)-\alpha(s))Q_sV(s))_s+(Q_s(\pas V(s)),\Theta(s)Q_sV(s))_s\\
 & = & (Q_s(A(s)-\alpha(s))V(s),(A(s)-\alpha(s))Q_sV(s))_s-\\ & & -(\tilde{\alpha}(s)Q_sV(s),(A(s)-\alpha(s))Q_sV(s))_s+\\
 & & +(H(s),(A(s)-\alpha(s))Q_sV(s))_s+(Q_s(\pas V(s)),\Theta(s)Q_sV(s))_s\\
 & = & \|(A(s)-\alpha(s))Q_sV(s)\|^2_s+(Q_s(A(s)-\alpha(s))\pi_sV(s),(A(s)-\alpha(s))Q_sV(s))_s-\\
 & & -(\tilde{\alpha}(s)Q_sV(s),(A(s)-\alpha(s))Q_sV(s))_s+(H(s),(A(s)-\alpha(s))Q_sV(s))_s+\\
 & & +(Q_s(\pas V(s)),\Theta(s),Q_sV(s))_s\\
 & \ge & \|(A(s)-\alpha(s))Q_sV(s)\|^2_s-c\,e^{-|\lambda|s}\|(A(s)-\alpha(s))Q_sV(s)\|_s-\\
 & & -\ve(s)\|Q_sV(s)\|_s\|(A(s)-\alpha(s))Q_sV(s)\|_s-c\,e^{-|\lambda|s}\|Q_sV(s)\|_s.
\end{eqnarray*}
Using lemma \ref{5.1.12A.} we obtain now for large $s$:
\begin{eqnarray*}
g''(s) & \ge & \|(A(s)-\alpha(s))Q_sV(s)\|^2_s-c\,e^{-|\lambda|s}\|(A(s)-\alpha(s))Q_sV(s)\|_s-\\
& & -\ve(s)\|Q_sV(s)\|_s\|(A(s)-\alpha(s))Q_sV(s)\|_s-c\,e^{-|\lambda|s}\|Q_sV(s)\|_s-\\
 & & -\ve(s)\|Q_sV(s)\|^2_s\\
 & = & \|(A(s)-\alpha(s))Q_sV(s)\|_s\Big(\|(A(s)-\alpha(s))Q_sV(s)\|_s-\ve(s)\|Q_sV(s)\|_s\Big)-\\
 & & -c\,e^{-|\lambda|s}\|(A(s)-\alpha(s))Q_sV(s)\|_s-c\,e^{-|\lambda|s}\|Q_sV(s)\|_s-\\
 & & -\ve(s)\|Q_sV(s)\|^2_s\\
 & \ge & (\delta-\ve(s))\|(A(s)-\alpha(s))Q_sV(s)\|_s\|Q_sV(s)\|_s-\\
 & & -c\,e^{-|\lambda|s}(1+\delta)\|Q_sV(s)\|_s-\ve(s)\|Q_sV(s)\|^2_s\\
 & \ge & \tilde{\delta}^2g(s)-c\,e^{-|\lambda|s}
\end{eqnarray*}
for a suitable positive number $\tilde{\delta}$. We remark that $g(s)$ converges to zero as $s\rightarrow\infty$ since the remainder in the asymptotic formula and all its derivatives do. We introduce now the function
\[
\beta(s):=g(s)+\frac{ce^{-|\lambda|s}}{|\lambda|^2-\tilde{\delta}^2},
\]
which also tends to zero as $s\rightarrow \infty$. We have
\[
\beta''(s)\ge \tilde{\delta}^2\beta(s).
\]
Defining $\gamma(s):=\beta(s)-\beta(s_0)e^{-\tilde{\delta}(s-s_0)}$ we get $\gamma''(s)\ge \tilde{\delta}^2\gamma(s)$, $\gamma(s)\rightarrow 0$ as $s\rightarrow\infty$ and $\gamma(s_0)=0$ which implies that $\gamma$ is a non positive function. Therefore,
\begin{equation}\label{estimate-for-the-function-g}
g(s)\le \beta(s_0)e^{-\tilde{\delta}(s-s_0)}+\frac{ce^{-|\lambda|s}}{||\lambda|^2-\tilde{\delta}^2|}\le c\,e^{-\hat{\delta}s}
\end{equation}
for suitable positive constants $c,\hat{\delta}$. We now have to show exponential decay for $|\pi_sV(s)|$ and all the derivatives of $\alpha(s)-\lambda$. The proof will be by induction with respect to $N$, the length of the vector $V(s)$. We start by establishing the desired estimates for the case $N=1$. We claim that 
\[
|(\pas\xi(s),\xi(s))_s|\le c\,e^{-|\lambda|s},
\]
which follows easily from $\|\xi(s)\|_s=1$ since 
\[
0=(\pas\xi(s),\xi(s))_s+(\xi(s),\Gamma_1(s)\xi(s))_s,
\]
where $\|\Gamma_1(s)\|$ has the above exponential decay. We conclude that $|(\pas\xi(s),\pi_s\xi(s))_s|$ also decays exponentially since $|(\pas\xi(s),Q_s\xi(s))_s|$ does. We calculate
\begin{eqnarray*}
(\pas\xi(s),\pi_s\xi(s))_s & = & \frac{(e,\xi(s))_s}{\|e\|^2_s}(\pas\xi(s),e)_s\\
 & = & \frac{(e,\xi(s))_s}{\|e\|^2_s}(A(s)\xi(s),e)_s-\alpha(s)\frac{(e,\xi(s))^2_s}{\|e\|^2_s}\\
 & = & \frac{(e,\xi(s))_s}{\|e\|^2_s}\Big((\xi(s),\Theta(s)e)_s+(\xi(s),(A(s)-A_{\infty})e)_s+\\
& & (\xi(s),\lambda e)_s-\alpha(s)(\xi(s),e)_s\Big)
\end{eqnarray*}
and
\begin{eqnarray*}
(\lambda-\alpha(s)) & = & \frac{\|e\|^2_s}{(e,\xi(s))_s^2}(\pas\xi(s),\pi_s\xi(s))_s+\\
 & & +\frac{1}{(\xi(s),e)_s}(\xi(s),(A_{\infty}-A(s))e-\Theta(s)e)_s.
\end{eqnarray*}
Recalling that $\inf_s(\xi(s),e)_s>0$ we conclude that for suitable constants $c,\hat{\delta}>0$
\begin{equation}\label{expo-decay-for-lambda-alpha}
|\lambda-\alpha(s)|\le c\,e^{-\hat{\delta}s}.
\end{equation}
We compute
\begin{eqnarray*}
\pi_s\pas\xi(s) & = & \pi_s(A(s)-\alpha(s))\xi(s)\\
 & = & \frac{1}{\|e\|_s^2}((A(s)-\alpha(s))\xi(s),e)_s\cdot e\\
& = & -\alpha(s)\frac{(\xi(s),e)_s}{\|e\|_s^2}\cdot e+\frac{1}{\|e\|_s^2}(\xi(s),A(s)e)_s\cdot e+\\
 & & +\frac{1}{\|e\|_s^2}(\xi(s),\Theta(s)e)_s\cdot e\\
 & = & -\alpha(s)\frac{(\xi(s),e)_s}{\|e\|_s^2}\cdot e+\frac{1}{\|e\|_s^2}(\xi(s),(A(s)-A_{\infty})e)_s\cdot e+\\
 & & +\frac{1}{\|e\|_s^2}(\xi(s),\lambda e)_s\cdot e+\frac{1}{\|e\|_s^2}(\xi(s),\Theta(s)e)_s\cdot e
\end{eqnarray*}
so that with (\ref{expo-decay-for-lambda-alpha}) and (\ref{estimate-for-the-function-g})
\begin{eqnarray}\label{expo-decay-for-pi_s-xi}
\|\pas\xi(s)\|_{L^2([0,1])} & \le & \|Q_s(\pas\xi(s))\|_{L^2([0,1])}+\|\pi_s(\pas\xi(s))\|_{L^2([0,1])} \nonumber \\
 & \le & c\,e^{-\hat{\delta}s}.
\end{eqnarray} 
Now
\[
\pi_s(\xi(s)-e)=\frac{(\xi(s),e)_s}{\|e\|_s^2}\cdot e-e=-\int_s^{\infty}\frac{d}{d\sigma}\frac{(\xi(\sigma),e)_{\sigma}}{\|e\|_{\sigma}^2}\cdot e\,d\sigma,
\]
and the integrand has exponential decay by our previous estimates. Since we have already shown that $\|Q_s\xi(s)\|_s$ decays exponentially, we obtain
\begin{equation}\label{expo-decay-of-xi-e}
\|\xi(s)-e\|_{L^2([0,1])}\le c\,e^{-\hat{\delta}s}
\end{equation}
for suitable constants $c,\hat{\delta}>0$. We can now complete the proof by induction as follows: Differentiating equation (\ref{5.1.21.}) for $\alpha(s)$ we obtain
\begin{eqnarray*}
\alpha'(s) & = & (\pas\xi(s),\Gamma_1(s)\xi(s))_s+(\xi(s),\Gamma_1(s)\pas\xi(s))_s+(\xi(s),\Gamma^2_1(s)\xi(s))_s+\\
 & & +(\xi(s),\pas\Gamma_1(s)\xi(s))_s+2(\pas\xi(s),A(s)\xi(s))_s+(\pas\xi(s),\Theta(s)\xi(s))_s+\\
 & & +(\xi(s),\pas M(v)\cdot M(v)\,A(s)\xi(s))_s.
\end{eqnarray*}
All the terms containing $\Gamma_1,\Theta$ or derivatives of $M(v)$ already decay at an exponential rate and will continue to do so if differentiated. We will summarize all those by $H(s)$.
Substituting $A(s)\xi(s)=\pas\xi(s)+\alpha(s)\xi(s)$ and using $(\nabla_s\xi(s),\xi(s))_s=0$, we then obtain
\begin{eqnarray*}
\alpha'(s) & = & 2\|\pas\xi(s)\|^2_s+2\alpha(s)(\pas\xi(s),\xi(s))_s+H(s)\\
 & = & 2\|\pas\xi(s)\|^2+H(s).
\end{eqnarray*}
Hence exponential decay of all derivatives $\|\pas^l\xi(s)\|_{L^2([0,1])}$ up to order $k\ge 1$, implies exponential decay of the derivative $\alpha^{(k+1)}(s)$. Denoting exponentially decaying expressions by $H(s)$, the PDE for $V(s)$ yields
\begin{eqnarray*}
\pi_s(\pas V(s)) & = & \pi_s(A(s)-\alpha(s))V(s)-\tilde{\alpha}(s)\pi_sV(s)+H(s)+E(s)\\
 & = & \frac{1}{\|e\|_s^2}[(\lambda-\alpha(s))(V(s),e)_s\, e+(V(s),(A(s)-A_{\infty})e)_s\, e]-\\
 & & -\tilde{\alpha}(s)\pi_sV(s)+H(s)+E(s),
\end{eqnarray*}
i.e. exponential decay of $(\|\pas^k\xi(s)\|_{L^2([0,1])})_{0\le k\le N-1}$ and $\left(\frac{d^k}{ds^k}(\lambda-\alpha(s))\right)_{0\le k\le N-1}$ implies exponential decay of $\|\pi_s\pas^N\xi(s)\|_{L^2([0,1])}$ and therefore of $\|\pas^N\xi(s)\|_{L^2([0,1])}$ in view of (\ref{estimate-for-the-function-g}). By iteration we obtain exponential decay for all derivatives of $\alpha$ and the $L^2$--norms of all s--derivatives of $\xi(s,t)$. Using the PDE for $V$, (\ref{PDE-for-V}), we also obtain exponential decay of $\|\pas^k\pat^l\xi(s)\|_{L^2([0,1])}$, and the Sobolev embedding theorem finally implies the assertion of the theorem.

\qed\\

\subsection{The asymptotic formula in local coordinates}
We will express the asymptotic formula in theorem \ref{asymptotic-formula-theorem} in coordinates near $\{0\}\times{\mathcal L}$ for later reference.\\
Recall that we have used proposition \ref{3.1.1.} and the modification (\ref{make-boundary-conditions-flat}) to derive the following coordinates on suitable neighborhoods $V_{\pm}$ of the points $p_{\pm}\in{\mathcal L}$:
\begin{equation}\label{make-boundary-conditions-flat-2}
\psi_{\pm}:{\bf R}^4\supset B_{\eps}(0)\widetilde{\longrightarrow}V_{\pm}\subset\RM,
\end{equation}
\[
\psi_{\pm}(0)=p_{\pm},
\]
\[
\psi_{\pm}({\bf R}^2\times\{0\}\times\{0\})=(\RL)\cap V_{\pm},
\]
\[
\psi_{\pm}(\{0\}\times{\bf R}\times\{0\}\times{\bf R}^{\pm})=(\{0\}\times {\mathcal D})\cap V_{\pm}.
\]
Using the coordinates $(\tau,\theta,x,y)$ for ${\bf R}^4$, the contact form on $\{0\}\times{\bf R}^3$ is then given by
\[
\hat{\lambda}_{\pm}=\psi_{\pm}^{\ast}\lambda=dy+\left(x+q(\theta)y\right)d\theta\ ,\ q(\theta):=\frac{a(\theta)}{b(\theta)}
\]
with Reeb vector field
\[
X_{\hat{\lambda}_{\pm}}=\frac{\partial}{\partial y}-q(\theta)\frac{\partial}{\partial x}
\]
(recall that the functions $a,b$ determine how the surface ${\mathcal D}$ is wrapping itself around the knot ${\mathcal L}$, see proposition \ref{3.1.1.}). Then $v_{\pm}(s,t):=(\psi^{-1}_{\pm}\circ \tu_0)(s,t)$ is the representative of the pseudoholomorphic curve $\tu_0$ in the above coordinates near the ends and the results of this section \ref{asymptotic-behaviour} (exponential decay estimates, asymptotic formula etc.) all refer to the maps $v_{\pm}$. We will now compute the eigenvector $e(t)$ in theorem \ref{asymptotic-formula-theorem} explicitly in the above coordinates.\\

The vectors
\[
\hat{e}_1(\theta,x,y):=(0, 1, -q'(\theta)y+q(\theta)(x+q(\theta)y), -(x+q(\theta)y))
\]
and
\[
\hat{e}_2(\theta,x,y) := (0, 0, 1, 0)
\]
generate the contact structure $\ker\hat{\lambda}(\theta,x,y)$ so that the almost complex structure $\hat{J}$ induced by $\tilde{J}$ is given by
\begin{equation}\label{matrix-of-J-hat}
\hat{J}(\tau,\theta,x,y)\hat{e}_1(\theta,x,y)=-\hat{e}_2(\theta,x,y)\ ,\ \hat{J}(\tau,\theta,x,y)
\hat{e}_2(\theta,x,y)=\hat{e}_1(\theta,x,y),
\end{equation}
\[
\hat{J}(\tau,\theta,x,y)(1,0,0,0)=(0,0,-q(\theta),1)\ ,\ \hat{J}(\tau,\theta,x,y)(0,0,-q(\theta),1)=(-1,0,0,0)
\]
or
\[
\hat{J}(\tau,\theta,x,y) =  \left(\begin{array}{cc}
0 & -(x+q(\theta)y) \\ 0 & yq'(\theta) \\
-q(\theta) & -1+yq'(\theta)((x+q(\theta)y)q(\theta)-yq'(\theta)) \\1 & -(x+q(\theta)y)yq'(\theta)
\end{array}\right. 
\]
\[
\left.\begin{array}{cc}
0 & -1 \\1 & q(\theta) \\(x+q(\theta)y)q(\theta)-yq'(\theta) & q(\theta)((x+q(\theta)y)q(\theta)-yq'(\theta)) \\
-(x+q(\theta)y) & -(x+q(\theta)y)q(\theta)
\end{array}\right).
\]
Theorem \ref{asymptotic-formula-theorem} then provides the following formula for the map $v_{\pm}(s,t)$ if $|s|$ is large:
\[
v_{\pm}(s,t)=e^{\int_{s_0}^s\alpha_{\pm}(\tau)d\tau}[e_{\pm}(t)+r_{\pm}(s,t)],
\]
where $e_{\pm}(t)$ is an eigenvector of the operator
\[
A_{\pm\infty}:L^2([0,1],{\bf R}^4)\supset H^{1,2}_{L}([0,1],{\bf R}^4)\longrightarrow L^2([0,1],{\bf R}^4)
\]
\[
(A_{\pm\infty}\gamma)(t) = -\hat{J}(0,0,0,0)\dot{\gamma}(t)
\]
corresponding to some eigenvalue $\lambda_{\pm}$, and we saw earlier in proposition \ref{spectrum-of-A-infty} that $\lambda_{\pm}$ are integer multiples of $\pi/2$. In fact, we will mostly be concerned with the case where $\lambda_+=-\frac{\pi}{2}$ and $\lambda_-=\frac{\pi}{2}$. The subscript $'L'$ refers to the boundary conditions $\gamma(0)\in {\bf R}^2\times\{0\}\times\{0\}$ and $\gamma(1)\in \{0\}\times{\bf R}\times\{0\}\times{\bf R}$. The matrix $\hat{J}(0,0,0,0)$ is given by 
\[
\hat{J}(0,0,0,0)=\left(\begin{array}{cccc} 0 & 0 & 0 & -1 \\ 0 & 0 & 1 & q(0) \\ -q(0) & -1 & 0 & 0 \\ 1 & 0 & 0 & 0 \end{array}\right)
\]
so that we have to solve the following system of differential equations for $e_{\pm}=(e_1,\ldots,e_4):[0,1]\rightarrow {\bf R}^4$
\begin{eqnarray*}
\dot{e}_1(t) & = & -\lambda_{\pm} e_4(t)\\
\dot{e}_2(t) & = & \lambda_{\pm}(e_3(t)+q(0)e_4(t))\\
\dot{e}_3(t) & = & -\lambda_{\pm}(q(0)e_1(t)+e_2(t))\\
\dot{e}_4(t) & = & \lambda_{\pm} e_1(t)
\end{eqnarray*}
with the boundary condition
\[
e_3(0)=e_4(0)=0\ ,\ e_1(1)=e_3(1)=0.
\]
If $\lambda$ is an integer multiple of $\pi$, we have 
\begin{equation}\label{even-eigenvector}
e(t)=\kappa(0, \cos(\lambda t), -\sin(\lambda t), 0)\, ,\, \kappa\neq 0
\end{equation}
Otherwise, if $\lambda_{\pm}$ is an odd integer multiple of $\pi/2$ then
\begin{equation}\label{odd-eigenvector}
e_{\pm}(t)=-{\kappa}_{\pm}(\cos(\lambda_{\pm} t), -q_{\pm}(0)\cos(\lambda_{\pm} t), 0, \sin(\lambda_{\pm} t))
\end{equation}
for some constants ${\kappa}_{\pm}\neq 0$. The asymptotic formula of theorem \ref{asymptotic-formula-theorem} then looks as follows:
\begin{eqnarray}\label{explicit-asymptotic-formula}
v_{\pm}(s,t) & = & -{\kappa}_{\pm}e^{\int_{s_0}^s\alpha_{\pm}(\tau)d\tau}\Big(\cos(\lambda_{\pm} t), -q_{\pm}(0)\cos(\lambda_{\pm} t), 0, \sin(\lambda_{\pm} t)\Big)+ \nonumber \\
 & & +e^{\int_{s_0}^s\alpha_{\pm}(\tau)d\tau}\ve_{\pm}(s,t).
\end{eqnarray}
In the following we will denote by $\varepsilon(s,t)$ any ${\bf R}^4$-- or real--valued function which converges to zero with all its derivatives uniformly in $t$ as $s\rightarrow\pm\infty$ if we are not interested in the particular function. In order to simplify notation we will often drop the subscript $\pm$. Using proposition \ref{alphaprime-goes-to-zero} we obtain the following asymptotic formulas for the derivatives of $v(s,t)$
\begin{eqnarray}\label{explicit-asymptotic-formula-for-us}
\partial_s v(s,t) & = & e^{\int_{s_0}^s\alpha(\tau)d\tau}\cdot \nonumber\\
 & & \cdot\Big[-{\kappa}\Big(\lambda\cos(\lambda t), -\lambda q(0)\cos(\lambda t), 0, \lambda\sin(\lambda t)\Big)+\varepsilon(s,t)\Big],
\end{eqnarray}
\begin{eqnarray}\label{explicit-asymptotic-formula-for-ut}
\partial_t v(s,t) & = & e^{\int_{s_0}^s\alpha(\tau)d\tau}\cdot \nonumber\\
 & & \cdot\Big[-{\kappa}\Big(-\lambda\sin(\lambda t), \lambda q(0)\sin(\lambda t), 0, \lambda\cos(\lambda t)\Big)+\varepsilon(s,t)\Big].
\end{eqnarray}

We will sometimes use the coordinates given by proposition \ref{3.1.1.} without making the boundary conditions 'flat' as in (\ref{make-boundary-conditions-flat}). In this case the appropriate versions of (\ref{explicit-asymptotic-formula}) and (\ref{explicit-asymptotic-formula-for-us}) are the following. If $\lambda_{\pm}$ is an odd integer multiple of $\pi/2$ we have:
\begin{eqnarray}\label{odd-explicit-asymptotic-formula-v2}
 & & v_{\pm}(s,t)\nonumber \\ 
& = & -{\kappa}_{\pm}e^{\int_{s_0}^s\alpha_{\pm}(\tau)d\tau}\Big(\cos(\lambda_{\pm} t), -q_{\pm}(0)\cos(\lambda_{\pm} t), q_{\pm}(0)\sin(\lambda_{\pm} t), \sin(\lambda_{\pm} t)\Big)+ \nonumber \\
 & & +e^{\int_{s_0}^s\alpha_{\pm}(\tau)d\tau}\ve_{\pm}(s,t)
\end{eqnarray}
and
\begin{eqnarray}\label{odd-explicit-asymptotic-formula-for-us-v2}
& & \partial_s v_{\pm}(s,t)\nonumber \\
 & = & e^{\int_{s_0}^s\alpha_{\pm}(\tau)d\tau}\cdot \Big[-{\kappa}_{\pm}\Big(\lambda_{\pm}\cos(\lambda_{\pm} t), -\lambda_{\pm} q_{\pm}(0)\cos(\lambda_{\pm} t),\\
 & &  \lambda_{\pm} q_{\pm}(0)\sin(\lambda_{\pm} t), \lambda_{\pm}\sin(\lambda_{\pm} t)\Big)+\varepsilon_{\pm}(s,t)\Big].\nonumber
\end{eqnarray}
For $\lambda_{\pm}\in{\bf Z}\pi$ we have
\begin{eqnarray}\label{even-explicit-asymptotic-formula-v2}
 & & v_{\pm}(s,t)\nonumber \\ 
& = & {\kappa}_{\pm}e^{\int_{s_0}^s\alpha_{\pm}(\tau)d\tau}\Big(0 , \cos(\lambda_{\pm} t), -\sin(\lambda_{\pm} t), 0\Big)+ \nonumber \\
 & & +e^{\int_{s_0}^s\alpha_{\pm}(\tau)d\tau}\ve_{\pm}(s,t)
\end{eqnarray}
and
\begin{eqnarray}\label{even-explicit-asymptotic-formula-for-us-v2}
\partial_s v_{\pm}(s,t) & = & e^{\int_{s_0}^s\alpha_{\pm}(\tau)d\tau}\cdot \nonumber\\
 & & \cdot\Big[{\kappa}_{\pm}\Big(0 , \lambda_{\pm} \cos(\lambda_{\pm} t), -\lambda_{\pm}\sin(\lambda_{\pm} t),0\Big)+\varepsilon_{\pm}(s,t)\Big].
\end{eqnarray}

\end{document}